\RequirePackage[leqno]{amsmath}
\documentclass[leqno]{amsart}
\makeatletter
\usepackage[utf8]{inputenc}     % Accept different input encodings
\usepackage[T1]{fontenc}
\usepackage{multicol}
\usepackage{enumitem}
\usepackage{color, graphics}
\usepackage[demo]{graphicx}
\usepackage{efbox,graphicx}
\efboxsetup{linecolor=gray,linewidth=0.8pt}
\usepackage{floatrow}
\usepackage{pgfplots}
\usepackage{graphicx,multicol,latexsym,amsmath,amssymb}
\usepackage{pgfplots}
\usepackage{tikz} 
\usepackage{amsfonts}
\usepackage{enumitem}
\setlist[itemize]{topsep=0pt,after=\vspace{1.5\baselineskip}}
 \usepackage[colorlinks=true]{hyperref}
\usepackage{empheq}
\usepackage{tabularx}
\usepackage[leftcaption]{sidecap}
\sidecaptionvpos{figure}{m}
\usepackage{tikz}
\usepackage[margin=1in]{geometry}
\usepackage{float}
\usepackage{subcaption}
\usetikzlibrary{backgrounds,automata}

\newcommand{\into}{\int_\Omega}

\newcount\rc@count
\let\rc@clearconstantlist\empty
\newcommand\rc@clearconstant[1]{\global\expandafter\let\csname rc@const@#1\endcsname\undefined}
\newcommand\resetconstants[1]{%
    \def\rc@constname{#1}% Set the new base name of the constants to the argument
    \global\rc@count=1\relax % Reset the constant counter to 1
    \bgroup 
        \let\\\rc@clearconstant % map over the list of constants that have been defined, clearing each of them.
        \rc@clearconstantlist
        \global\let\rc@clearconstantlist\empty % Globally empty the list of constants.
    \egroup
}
\newcommand\const[1]{%
    \@ifundefined{rc@const@#1}{%
        \expandafter\xdef\csname rc@const@#1\endcsname{%
           \noexpand\rc@useconst{\rc@constname}{\the\rc@count}%
        }%
        \g@addto@macro\rc@clearconstantlist{\\{\mathrm{#1}}}%
        \global\advance\rc@count1\relax
    }{}%
    \csname rc@const@#1\endcsname
}
\DeclarePairedDelimiter\abs{\lvert}{\rvert}

\DeclarePairedDelimiter\tonda{(}{)}

\newcommand\rc@useconst[2]{{#1}\textsubscript{#2}}
\setlist[itemize]{noitemsep, topsep=0pt}

\def\R{\mathbb R} \def\N{\mathbb N}

\def\R{\mathbb R} \def\N{\mathbb N} 
\def\TM{T_{\rm{max}}} 
\def
\@cite
#1#2{[#1\if@tempswa, #2\fi]}
\makeatother
\addcontentsline{toc}{section}{Abstract}
\newtheorem{theorem}{Theorem}[section]
\newtheorem{corollary}[theorem]{Corollary}

\newtheorem{lemma}[theorem]{Lemma}

\newtheorem{remark}{Remark}
\newtheorem{convention}{Convention}

\newcounter{cnstcnt}

\title[Boundedness in chemotaxis models with different dissipative actions]{Dissipation through combinations of nonlocal and gradient nonlinearities
in chemotaxis models}
\author[Rafael D\'iaz Fuentes, Silvia Frassu and Giuseppe Viglialoro]{}
\keywords{Chemotaxis, Global existence, Nonlocal terms, Gradient terms, Boundedness. \\
\textit{$^*$Corresponding author}: giuseppe.viglialoro@unica.it}

\begin{document}

\maketitle
\centerline{\scshape{\scshape{Rafael D\'iaz Fuentes, Silvia Frassu \and Giuseppe Viglialoro$^{*}$}}}
\medskip
{
\medskip
\centerline{$^\sharp$Dipartimento di Matematica e Informatica}
\centerline{Universit\`{a} di Cagliari}
\centerline{Via Ospedale 72, 09124. Cagliari (Italy)}
\medskip
}
\bigskip
\tableofcontents
\begin{abstract}
This work concerns with a class of chemotaxis models in which external sources, comprising nonlocal and gradient-dependent damping reactions, influence the motion of a cell density attracted by a chemical signal. The mechanism of the two densities is studied in bounded and impenetrable regions. In particular, it is seen that no gathering effect for the cells can appear in time provided that the damping impacts are sufficiently strong. 

Mathematically, we study this problem 
\begin{equation}\label{problem_abstract}
\tag{$\Diamond$}
\begin{cases}
u_t=\nabla \cdot \left((u+1)^{m_1-1}\nabla u -\chi u(u+1)^{m_2-1}\nabla v\right)+B(x,t,u,\nabla u)&{\rm in}\ \Omega \times (0, \TM),\\
\tau v_t=\Delta v-v+f(u) &{\rm in}\ \Omega \times (0, \TM),\\
u_\nu=v_\nu=0 &{\rm on}\ \partial\Omega \times (0, \TM),\\
u(x, 0)=u_0(x), \tau v(x,0)= \tau v_0(x) &x \in \bar{\Omega},
\end{cases} 
\end{equation}
for
\begin{equation*}%\label{ExpressionofBForSources}
B(x,t,u,\nabla u)=B=
\begin{cases}
au^\alpha-b u^\beta-c\int_\Omega u^\delta, \\
au^\alpha-b u^\alpha \int_\Omega u^\beta-c|\nabla u|^\delta,
\end{cases} 
\end{equation*}
and where $\Omega$ is a bounded and smooth domain of $\R^n$ ($n\in \N$), $\tau\in\{0,1\}$, $m_1,m_2\in \R$, $\chi,a,b>0$, $c\geq 0$, and $\alpha, \beta,\delta\geq 1$. Herein $u$ stands for the population density,  $v$ for the chemical signal, $f$ for a regular function describing the production law and $\TM$ for the maximal time of existence of any nonnegative classical solution $(u,v)$ to system \eqref{problem_abstract}, emanating from nonnegative and sufficiently regular initial distribution $u_0, v_0$. For each of the expressions of $B$, sufficient conditions on parameters of the models ensuring that 
%\begin{itemize}
%\item (the subquadratic case) $1\leq \alpha <2 \quad \textrm{and} \quad \beta>\frac{n+4}{2}-\alpha,$
%\item (the superquadratic case) $\beta>\frac{n}{2} \quad \textrm{and} \quad 2\leq \alpha < 1+ %\frac{2\beta}{n},$
%\end{itemize}
actually $\TM=\infty$ and $u$ and $v$ are uniformly bounded are established.

In the literature, most of the results concerning chemotaxis models with external sources deal with classical logistics, for which $B=a u^\alpha-b u^\beta$. Thereafter, the introduction of dissipative effects as those expressed in $B$ is the main novelty of this investigation. On the other hand, this paper extends the analyses in \cite{chiyoduzgunfrassuviglialoro2024,BianEtAlNonlocal,Latos2020}.
\end{abstract}
\resetconstants{c}
%\tableofcontents
\section{Introduction and motivations. Nonlocal and gradient-type sources}\label{IntroSection}
In this paper we will consider chemotaxis models perturbed by reactions involving dissipative effects associated to nonlocal and gradient nonlinearities. Subsequently, below we present some discussions concerning such sources in the framework of biological models formulated by Reaction/Diffusion and Taxis/Reaction/Di\-ffu\-sion equations.
\subsection{The Reaction/Diffusion equation in biology. Basic mathematical concepts}
Reaction/Di\-ffu\-sion equations (RDe, singular and plural) describe models that correspond to several physical phenomena. It is known that the diffusion term in RDe contributes to transfer the concentration of a physical quantity throughout the spatial habitat where it is distributed, whereas the reaction term (which can be of different nature) comes from the mass action laws. 

In mathematical biology, exactly accordingly to \cite{Murray1,Murray2}, for a 
particles' distribution $u= u(x,t)$, being $x$ a spatial position inside a habitat $\Omega$ occupied by the species itself and $t\in (0,\TM)$ a temporal instant, the evolution modeling its dynamic is a system of  RDe reading as
 \begin{equation}\label{GeneralRDe} 
 u_t=\nabla \cdot (D(u) \nabla u)+f(u)\quad \textrm{in } \Omega \times (0,\TM).
\end{equation}  
Herein, the diffusivity $D(u)$ influences the self-diffusion process through the term $\nabla \cdot (D(u) \nabla  u)$, whereas $f(u)$ represents the actions affecting the total amount of the species in the spatial domain and throughout the time. The above equation has to be equipped with nonnegative initial configuration $ u(x, 0) = u_0(x)$ and appropriate boundary conditions on $\partial \Omega$ (the boundary of $\Omega$), as for instance Dirichlet, no-flux (Neumann) and/or periodic. 

Classical questions connected to solutions to problem \eqref{GeneralRDe} concern their global existence (and, possibly, the boundedness) and their blow-up (at finite or infinite time). In particular, if $u=u(x,t)$ is a solution for any $x\in \Omega$ and $t\in (0,\TM)$ and $\TM=\infty$, then it is global (and bounded if $\sup_{t\in (0,\infty)}||u(\cdot,t)||_{L^\infty(\Omega)}< \infty$). If for some $\TM$ (finite or infinite) it occurs that $\limsup_{t\rightarrow \TM} ||u(\cdot,t)||_{L^\infty(\Omega)}=\infty$, then $u$ blows-up at $\TM$; if $\TM$ is finite, $\delta$-formations in the region appear.

%, etc. Different types of boundary conditions have been used in RD equations, depending
%on the nature of the problems; no-flux and periodic boundary conditions are among the most familiar
\subsection{The local and nonlocal Fisher--KPP equation as prototypes}\label{RDeOverviewSectionNonLocal} 
Nonlocalities play a crucial role in practical mathematical representations, especially at small scales: they define the attributes of two entities situated in separate locations. 

Motivated by the seminal work by Fisher on gene spread \cite{Fisher-1937-WAA}, researchers have developed various ways of modeling spatial interactions in biology.

In the specific, the Fisher--KPP equation in its original one spatial dimensional version, describes a mechanism modeling the growth of a single population as well as its wave propagation. It is a particular case of the above RDe \eqref{GeneralRDe}, for which $D(u)$ is constant, and given by
\begin{equation}\label{Fisher-KPPeq}
u_t- \Delta u = f(u), \quad \textrm{with } f(u)=u- u^\beta, \beta>1,
\end{equation} 
where $f$ can be interpreted as the competition between the growth rate of the population, i.e.,  $u$, and the death/damping one, i.e.,  $-u^\beta$. 

In its \textit{local} form \eqref{Fisher-KPPeq}, the population's segregation/consumption rates at any given position $x$ are proportional to the populations density itself at the same spatial point; in particular, the diffusion term contributes to the net density change in the spatial locations.

In \textit{nonlocal} phenomena, such rates at a spatial point $x$ may be influenced by factors in close locations (\cite{Furter-Grinfeld,BillinghamNonlocalFisher,Pal2024NonlocalMI}). In the specific, idealizing  the interaction between individuals located at the spatial point $y$ with those at position $x$ through the probability density function $\varphi(x-y)$, the local RDe \eqref{Fisher-KPPeq} is modified into this nonlocal integro-differential equation
\begin{equation*}
u_t- \Delta u = f(u,\varphi), \quad \textrm{with } f(u)=u- u* \varphi.
\end{equation*} 
Herein $*$ stands for the convolution operator, and the related term  captures the long-range interaction in the spatial domain, one of the key factors missing in the local model. 

From our mathematical perspective, nonlocal terms are expressed as an integral in space. These terms can also be combined with local terms, either in an additive or multiplicative manner, so to describe a competition between nonlocal sources and local
damping/dissipative actions. Specifically, for some proper $\alpha, \beta>0$ let us focus on this Fisher-KPP RDe
\begin{equation}\label{Fisher-KPPeqNonloc}
u_t=\Delta u + \int_\Omega u^\beta dx- u^\alpha \quad \textrm{in } \Omega \times (0,\TM),
\end{equation}
equipped with some initial distribution $u(x,0)=u_0(x)$ and some boundary conditions.  It is known from \cite{WangWang1996} that for homogeneous Dirichlet condition ($u=0$ on $\partial \Omega$) the value $\beta = \alpha$ represents a critical blow-up exponent, in the sense that for $\alpha>\beta$ (damping rate surpassing growth one) any initial data $u_0(x)$ emanates globally bounded solutions, whereas for $\beta>\alpha$ (growth rate surpassing damping one) there are initial distributions $u_0(x)$ leading to unbounded solutions. (See also \cite[Theorem 43.1]{QS} and \cite{Souplet-BlowUpNonlocal,Biler-P-Nonlocal,BianEtAl-FisherKPP} for further results in similar contexts.) 
%It appears worthwhile to underline that by inverting the action of the source $f(u)=\int_\Omega u^\beta dx- u^\alpha$ in problem \eqref{NonLocalEquation}, in the sense that the growth is associated to the local term and the consumption on the nonlocal, we are not aware of any blow-up results; we refer to \cite{} for boundedness questions in the case $\beta=\alpha$, under Neumann boundary conditions ($u_\nu=0$ on $\partial \Omega$, $\nu$ being the outward normal to $\partial \Omega$).
\subsection{Dissipative gradient terms in biological models}\label{RDeOverviewSectionGradient} 
In addition to nonlocal reactions as those discussed for equation \eqref{Fisher-KPPeqNonloc}, gradient-dependent sources may influence the dynamics of a RDe. As a matter of fact, according to \cite{Souplet_Gradient}, a biological species with density $u$ and occupying a certain habitat evolves in time by displacement, birth/reproduction and death. In particular, the births are described by a superlinear power of such a distribution, the natural deaths by a linear one and the accidental deaths by a function of its gradient. For some $\alpha,\delta>0$ this leads to the equation
\begin{equation}\label{EquationWithGradientTerm}
u_t=\Delta u +u^\alpha-u-|\nabla u|^\delta \quad \textrm{in } \Omega \times (0,\TM).
\end{equation}
%this leads (through the choice \textcolor{blue}{$F=F(u)=-\nabla u$} and $h=h(u,\nabla u)=u^\rho-u- |\nabla u|^\gamma$ in \eqref{ContinuityEq}) to $$, with $\rho,\gamma>1$.  
Similarly to what said for nonlocal sources, the study of gradient terms in parabolic equations, and their impact on the possibility of blow-up, is well known; 
for instance, for the Dirichlet problem associated to equation \eqref{EquationWithGradientTerm} in \cite{chipot_weissler,kawohl_peletier,fila} and \cite[Sect. IV]{QS} the interplay between the exponents $\alpha$ and $\delta$  preventing or implying  blow-up are established.  (We mention \cite{SoupletRecent2011,Souplet2005Book} as well, where surveys and more results on the topic can be found.)
\section{Chemotaxis models. An overview on Keller--Segel systems with different sources}
\subsection{The Taxis/Reaction/Diffusion equation in biology. The Keller--Segel model}
In some biological mechanisms, the dynamics of a species may also be affected by some convection action; in this case, the particles are transferred from a place to another of the habitat due to the presence of a transport vector. Mathematically, in the RDe \eqref{GeneralRDe} has to be introduced a vector $j=j(u)$, which transports the species $u$ along its direction; the final formulation results into this Taxis/Reaction/Diffusion equation (TRDe)
 \begin{equation}\label{GeneralTRDe} 
 u_t=\nabla \cdot (D(u) \nabla u-uj(u))+f(u)\quad \textrm{in } \Omega \times (0,\TM).
\end{equation}  
For the chemotaxis model introduced by Keller and Segel in 1970s (\cite{Keller-1971-MC,Keller-1971-TBC}), the vector $j=j(u)$ is proportional to the gradient of a chemical signal $v=v(x,t)$, distributed in $\Omega$; more precisely $j=j(u)=\chi \nabla v$, for some $\chi>0$. (Below we will mention cases where $j$ explicitly depends on $u$.) In such a situation, and in absence of the source $f$, the partial differential equation modeling the motion of  $u$ is
\begin{equation}\label{problemOriginalKS}
u_t=\Delta u-\chi  \nabla \cdot (u \nabla v)\quad  \textrm{in}\quad \Omega \times (0,\TM).
\end{equation}
It describes how a chemotactical impact of the (chemo)sensitivity ($\chi$) provided by the chemical signal $v$  may break the natural diffusion  (associated to the Laplacian operator, $\Delta u$) of the cells. Indeed, the term $-\nabla \cdot (u \chi \nabla v)$ models the transport of $u$ in the direction  $\chi \nabla v$, the negative sign indicating the attractive effect that $v$ has  on the cells (higher for $\chi$ larger and for an increasing amount of $v$). As a consequence, when $v$ is produced by the same cells, and in such a scenario $v$ obeys 
\begin{equation}\label{SecondEqKS}
v_t=\Delta v-v+u\quad \textrm{in}\quad \Omega \times (0,\TM),
\end{equation}
the attractive impact may be so efficient as to lead the cell density to a gathering process generating a \textit{chemotactic collapse (or blow-up)}. 
\subsection{An overview on the Keller--Segel model}
In contrast to the pure RDe, for which the reaction $f$ plays a pivotal role in the behaviour of the solutions (as previously outlined in the context of models \eqref{Fisher-KPPeqNonloc} and \eqref{EquationWithGradientTerm} in terms of $\alpha$ and $\beta$ and $\alpha$ and $\delta$, respectively), the situation concerning Taxis/Diffusion equations (and henceforth even in the absence of reactions) is more complex.  

Mathematically, it has been proved that solutions to the initial-boundary value problem associated with equations \eqref{problemOriginalKS} and \eqref{SecondEqKS} may be globally bounded in time or may blow-up at finite time. This depends on the mass (i.e., $\int_\Omega u_0(x)dx$) of the initial data, its specific configuration, and the value of the sensitivity $\chi$. In one-dimensional settings, all solutions are uniformly bounded in time. In contrast, for $n\ge 3$, it is possible to construct solutions that blow-up at finite time for any arbitrarily small mass $m=\int_\Omega u_0(x)dx>0$. Conversely, in the case $n=2$, the value $4\pi$ serves to differentiate between two scenarios. In the first instance, diffusion overcomes self-attraction (if $\chi m<4\pi$), while in the second, self-attraction dominates (if $\chi m>4\pi$). Consequently, all solutions are global in time, and initial data leading to assembling processes at finite time can be detected. A comprehensive overview of these analyses can be found in the following references: \cite{HerreroVelazquez,Nagai, nagai1995blow}. These works represent seminal contributions to the field. 
\subsection{An overview on the Keller--Segel model with logistic sources}
From the perspective of blow-up prevention (or boundedness enforcement, when explosion scenarios are not known), the introduction of different external agents is quite efficient in several (\textit{but not} in all) cases. More exactly, if the evolution of $u$ in equation \eqref{problemOriginalKS} is also influenced by the presence of logistic terms behaving as $au-bu^{\beta}$, for $a, b >0$ and $\beta>1$, mathematical intuition suggests that superlinear damping effects should benefit the boundedness of solutions. Actually,  the prevention of $\delta$-formations in the sense of finite time blow-up for 
\begin{equation*}%\label{KS}
u_t=\Delta u-\chi  \nabla \cdot (u \nabla v)+au-bu^\beta\quad  \textrm{in}\quad \Omega \times (0,\TM),
\end{equation*}
when coupled with some equation implying the segregation of $v$ with $u$ (for instance \eqref{SecondEqKS}), has been established only for large values of $b$  (if $\beta=2$, see \cite{TelloWinkParEl}, \cite{W0}), whereas for some value of $\beta$ next to $1$ a blow-up scenario was detected, first for dimension $5$ or higher \cite{WinDespiteLogistic}, (see also \cite{FuestCriticalNoDEA} for an improvement of \cite{WinDespiteLogistic}), but later also for $n\geq 3$, in  \cite{Winkler_ZAMP-FiniteTimeLowDimension}.
\subsection{An overview on the Keller--Segel model with nonlocal sources}
Likewise to the aforementioned external sources, impacts behaving as
\begin{equation}\label{NonLocalSources}
au^\alpha-bu^\alpha\int_\Omega u^\beta \quad a, b>0 \textrm{ and } \alpha,\beta\geq 1,
\end{equation}
model a competition between a birth contribution, favoring instabilities of the species (especially for large values of $a$), and a death one opportunely contrasting this instability (especially for large values of $b$). In this context,  
%naturally arise. 
%\begin{enumerate}[label=\textbf{$\mathcal{Q}$:},ref=$\mathcal{Q}$]
%\item a
%\item \label{Question}
in a biological mechanism governed by the equation
\begin{equation*}%\label{KS-NonLocal}
u_t=\Delta u-\chi  \nabla \cdot (u \nabla v)+au^\alpha-bu^\alpha\int_\Omega u^\beta \quad  \textrm{in}\quad \Omega \times (0,\TM),
\end{equation*}
the external damping source may suffice to enforce boundedness of solutions, even for any large initial distribution $u_0$, arbitrarily small $b>0$ and in any large dimension $n$.
% Are, conversely, some restrictions on $n$ and/or $a, b$, $\alpha, \beta, u_0$ required? 
%\item c. 
%\end{enumerate}
%To our knowledge, most of the analyses connected to the aforementioned questions can be found in the literature when the equation for $v$, describing how $u$ secretes $v$, is of elliptic type, i.e. for some $\gamma\geq 1$
%\begin{equation}\label{EquationEllipticNonlinearProd}
%0=\Delta v-v +u^\gamma\quad  \textrm{in}\quad \Omega \times (0,\TM). 
%\end{equation}
With the aim of giving an overview in this direction and in the framework of (zero-flux) Keller--Segel systems, let us consider for $\tau\in\{0,1\}$, and proper $m_1,m_2,\gamma\geq 1$ and $h=h(x,t)$, the initial-boundary value problem associated to this model
\begin{equation}\label{Tello}
\begin{cases}
u_t=\nabla \cdot \left((u+1)^{m_1-1}\nabla u-\nabla \cdot (\chi u (u+1)^{m_2-2}\nabla v\right)+f(u) &\quad  \textrm{in}\quad \Omega \times (0,\TM),\\
\tau v_t=\Delta v-v+u^{\gamma}+h &\quad  \textrm{in}\quad \Omega \times (0,\TM).\\
%u_\nu=v_\nu=0 &{\rm on}\ \partial\Omega,\\
%u(x, 0)=u_0(x), v(x,0)=v_0(x), &x \in \Omega,
\end{cases}
\end{equation}
Moreover, let the nonlocal term be 
\begin{equation}\label{NonLocalTermTello}
f(u):=u\left(a_0-a_1u^\alpha+a_2\int_\Omega u^\alpha dx\right),
\end{equation}
where $\alpha \geq 1$, $a_0,a_1>0$ and $a_2\in \R$. When the equation for the chemical $v$ is elliptic (i.e., $\tau=0$), a higher number of results compared to the situation where $\tau=1$ is available. Precisely, dealing with uniform-in-time boundedness of classical solutions, for problem \eqref{Tello}, with  $h\equiv 0$ and reaction $f$ as in 
\eqref{NonLocalSources} we can mention: 
%leads to 
\begin{itemize}
\item[$\triangleright$] for the special case where $m_1=\gamma=a=b=1$ and $m_2=2$ in \cite{BianEtAlNonlocal}, whenever these assumptions (with $\alpha\geq 1, \beta>1$) $n\geq3,$ $2\leq \alpha <
1+\frac{2\beta}{n}$ or $\frac{n+4}{2}-\beta<\alpha<2$ are complied; 
\item[$\triangleright$] in \cite{Latos2020} for the case  $m_1=a=b=1$, $\gamma\geq 1$ and $m_2 \geq 2$ tied by $\gamma+m_2-1\leq \alpha <
1+\frac{2\beta}{n}$ or $\frac{n+4}{2}-\beta<\alpha<\gamma+m_2-1$;
\item[$\triangleright$] for general choices of the parameters $m_1>0,m_2\geq 1,a=b>0$, for $\gamma=1$, under the hypotheses that $m_2+\frac{n}{2}(m_2-m_1)-\beta< \alpha <m_1+\frac{2}{n}\beta$ or  $\alpha=m_2+\frac{n}{2}(m_2-m_1)-\beta$ together with $b$ large enough (see \cite{TaoFang-NonlinearNonlocal}).
\end{itemize}
Another indication demonstrating the depth of the study within the context of models with stationary equations for the stimulus can be found in \cite{BianEtAlWholeSpaceNonLocal,ChenWangDocMat,Li-Viglialoro-DIE-Nonlocal},
investigating nonlocal problems in the whole space $\R^n$, which are analogous to those presented in model \eqref{Tello}.

Conversely, when the mechanism for the cells and the stimulus are both evolutive ($\tau=1$), we are only aware of \cite{Tello2021} and \cite{chiyoduzgunfrassuviglialoro2024}. In \cite{Tello2021}, the authors consider model \eqref{Tello} with source as in \eqref{NonLocalTermTello}, and globality and convergence to the steady state are established whenever the coefficients of the system satisfy
\begin{equation*}%\label{ConditionsTello}
\alpha+1>m_2-1+\gamma\textrm{ and } a_1-a_2|\Omega|>0.
\end{equation*}
 On the other hand, in the very recent contribution \cite{chiyoduzgunfrassuviglialoro2024}, model \eqref{Tello} is discussed for $m_1=\gamma=1$, $m_2=2$, $h\equiv 0$ and $f(u)$ as in \eqref{NonLocalSources}; 
%\begin{equation}\label{NonLocalTermChiyoEtAl}
%f(u):=au^\alpha-bu^\alpha\int_\Omega u^\beta, \quad a,b,\alpha,\beta>0;
%\end{equation}
boundendess of related classical solutions is achieved for (the subquadratic case) $1\leq \alpha <2 \quad \textrm{and} \quad \beta>\frac{n+4}{2}-\alpha,$ and for (the superquadratic case) $\beta>\frac{n}{2} \quad \textrm{and} \quad 2\leq \alpha < 1+ \frac{2\beta}{n}$. 
%\end{itemize}
%actually $\TM=\infty$ and $u$ and $v$ are uniformly bounded.
\subsection{An overview on the Keller--Segel model with gradient-dependent sources}
In contrast to chemotaxis phenomena influenced by nonlocal sources, the incorporation of gradient nonlinearities in such mechanisms represents a relatively novel trend, and the overall analysis is still in its infancy. In this regard, we can only cite  \cite{ViglialoroDifferentialIntegralEquations}, \cite{IshidaLankeitVigliloro-Gradient} and \cite{LiAcostaColumbuViglialoro}, where also real biological interpretations for the gradient terms are given. Specifically, by introducing the source $f(u,\nabla u)=a u^\alpha-b u^\beta-c|\nabla u|^\delta$ to describe a reaction with dissipative gradient terms, for the TRDe
\begin{equation}\label{OnlyATTractionproblemBis}
%\begin{cases}
u_t= \nabla \cdot ((u+1)^{m_1-1}\nabla u)-\chi \nabla  \cdot (u(u+1)^{m_2-1} \nabla v) +f(u,\nabla u)\quad \textrm{ in } \Omega \times (0,T_{\textrm{max}}),
%\end{cases}
\end{equation}
coupled with an equation for the chemical $v$ which is produced by $u$, lower bounds of the maximal existence time of given solutions are derived in two- and three-dimensional settings, when for $m_1=m_2=1$ some technical restrictions on $\alpha, \beta$ and $\delta$ are considered (\cite{ViglialoroDifferentialIntegralEquations}). In parallel, in the same context, \cite{IshidaLankeitVigliloro-Gradient} deals with the existence and boundedness of classical solutions to problem \eqref{OnlyATTractionproblemBis} with $m_1=m_2=1$, as well; in this sense, such boundedness of classical solutions is achieved in any dimension $n$, for any initial data, whenever $\beta>\alpha$ and $\frac{2n}{n+1}<\delta\leq 2.$ Finally, 
\cite{LiAcostaColumbuViglialoro} focuses on this more general TRDe
\begin{equation}\label{OnlyATTractionproblemTRIS}
%\begin{cases}
u_t= \nabla \cdot ((u+1)^{m_1-1}\nabla u-\chi u(u+1)^{m_2-1} \nabla v+\xi u(u+1)^{m_3-1} \nabla w) +f(u,\nabla u)\quad \textrm{ in } \Omega \times (0,T_{\textrm{max}}).
%\end{cases}
\end{equation}
In this case, the aforementioned flux $j(u)=-\chi u(u+1)^{m_2-1} \nabla v+\xi u(u+1)^{m_3-1} \nabla w$ counts with an attractive contribution, associated to the chemottractant $v$ and the repulsive counterpart, connected to a further chemical $w$, called chemorepellent. 
Essentially \cite{LiAcostaColumbuViglialoro} generalizes the result in \cite{IshidaLankeitVigliloro-Gradient}, where the linear and only attraction version (i.e., $m_1=m_2=1$, $\xi=0$) of model \eqref{OnlyATTractionproblemTRIS} is addressed.
\subsection{Gathering processes in chemotaxis models with nonlocal and gradient-dependent nonlinearties}
For RDe with nonlocal and/or gradient nonlinearities, in $\S$\ref{RDeOverviewSectionNonLocal} and $\S$\ref{RDeOverviewSectionGradient} blow-up situations have been discussed. As to this kind of instabilities in TDRe in the framework of chemotaxis systems with nonlocal terms, we refer to  \cite{DuLiu-BlowUpNonlocalChemo,KavallarisEtAlBlowUpChemoNonlocal,WolanskyBlowUpNonlocalChemo}. Conversely, when gradient-dependent sources are involved in chemotaxis models we are not aware of blow-up scenarios; despite that, there are numerical evidences indicating that explosions may appear; \cite[$\S$7]{LiAcostaColumbuViglialoro}.

%Per le equazioni di reazione diffusione con termine non locale a con gradient terms, nella sezione 3 abbiamo citato situazioni dove apparivano fenomeni di esplosione. A questi lavoro potremmo aggiungere altri (a parte quella lista sui nonlocal cercare altro su gradient): ....Non ci sono invece risultati che confermino che la stessa situazione puo verificarsi nel caso di equazioni tipo chemotassi. Al contrario, alcune evidenze numeriche dimostrano che anche in questo caso puo esserci blow-up.

%Additionally, the suppression of some of the conditions in \eqref{ConditionsTello}, might provide (at least from the numerical point of view) some blow-up solution.
Such indications supporting blow-up mechanisms even for damping nonlocal/gradient-dependent terms, make meaningful the analysis concerning the suppression of these coalescence phenomena. In this sense, as specified below, we here derive conditions on the parameters defining a class of chemotaxis models perturbed by damping sources ensuring boundedness of related classical solutions. 
%We start with local RD systems applied to ecological and chemical interactions. In particular,
%how the nonlocal Fisher-KPP equation comes in a single species population is discussed. Then, the nonlocal
%models capture different kinds of nonlocal phenomena. 
%\begin{equation}\label{problemFisherKKP}
%\begin{cases}
%u_t=\Delta u+ u^\alpha-u^\alpha\int_\Omega u &{\rm in}\ \Omega \times (0, \TM),\\
%%v_t=\Delta v-v+u &{\rm in}\ \Omega \times (0, \TM),\\
%u_\nu=0&{\rm on}\ \partial\Omega \times (0, \TM),\\
%u(x, 0)=u_0(x) &x \in \bar{\Omega},
%\end{cases}
%\end{equation}
%so directly obtainable for $\chi=0$ and $\beta=a=b=1$ from problem \eqref{problem}. It is established that for bounded and smooth domains $\Omega$ of $\R^n$, with $n\in \N$, whenever $u_0=u_0(x)$ is sufficiently regular and such that $\int_\Omega u_0(x)dx<1$, then system \eqref{problemFisherKKP} admits global solutions for $n = 1, 2$ with any $1 \leq \alpha <2$,  or $n \geq 3$ with any $1 \leq  \alpha < 1 + 2/n.$ 
%\subsection{Chemotaxis models witn nonlinear diffusion and sensitivity}

%%%%%%%%%%%%%
\section{Description and contextualization of the investigated chemotaxis models. Main results}
\subsection{The models and their context in the literature}
The literature concerning chemotaxis models abounds enormously, and a detailed description does not seem appropriate, at least in this context: despite that, we can say that the more common variants of model \eqref{problemOriginalKS} are connected to different choices of $D$ and $j$ in equation \eqref{GeneralTRDe}, in presence or absence of $f$.
%logistic type sources with dampening effects on the cells' increasing (\cite{Lankeit,LankeitWangConsumptLogistic,TelloWinkParEl,ViglialoroBoundnessVeryWeak}), involve 
In the specific, when coupled with equations as (or similar to) the second in  \eqref{Tello}, different results incorporating nonlinear expressions of the diffusion $D(u)\simeq u^{m_1}$ and sensitivity $j(u)\simeq u^{m_2}\nabla v$, are available: see, for instance,  \cite{CieslakStinnerFiniteTime,CieslakStinnerNewCritical,HashiraIshidaYokotaJDEBlow-UP,YokotaEtAlNonCONVEX,TaoWinkParaPara}. Additionally, also the degree of knowledge for nonlinear chemotaxis models with logistics $f$ as those described in \eqref{Fisher-KPPeq} is rather rich: we mention, as examples, \cite{GalakhovTelloJDE-2016,WangZhengJDE-2014,HuTao-2017}, and references therein. 

Conversely, \textit{this article aims at putting the emphasis on a much less in-depth scenario, wanting in particular to analyze the impact of nonlocal and gradient-type terms on the dynamics of specific chemotaxis models.} 
%, consider weaker (stronger) laws for the production of the chemical signal (\cite{LankeitViglialoroAAM,LiuTaoFullyParNonlinearProd,WinklerNoNLinearanalysisSublinearProduction}) or, more generally, suitable combinations of (some of) these actions (\cite{CaoZheng,ViglialoroWoolley2018,ZhengJJDE}).  
%Moreover, for the same boundedness purpose, other models take into account an additional chemical signal, produced or absorbed by the same cells but repelling them, thus influencing the overall dynamics: in this way, as explained at the beginning of the section, model \eqref{problem} is encompassed in the above casuistry.
%In this paper we consider 
%\begin{equation}\label{problem}
%\begin{cases}
%u_t=\Delta u-\chi\nabla \cdot (u\nabla v)+au^\alpha-bu^\alpha\int_\Omega u^\beta &{\rm in}\ \Omega \times (0, \TM),\\
%v_t=\Delta v-v+u &{\rm in}\ \Omega \times (0, \TM),\\
%u_\nu=v_\nu=0 &{\rm on}\ \partial\Omega \times (0, \TM),\\
%u(x, 0)=u_0(x), v(x,0)=v_0(x) &x \in \bar{\Omega},
%\end{cases}
%\end{equation}

As a matter of fact, the main interest of our research is to investigate these models:
\begin{equation}\label{problem_int}
\tag{$\mathcal{A}$}
\begin{cases}
u_t=\nabla \cdot \left((u+1)^{m_1-1}\nabla u -\chi u(u+1)^{m_2-1}\nabla v\right)+a u^\alpha-b u^\beta-c \int_\Omega u^\delta &{\rm in}\ \Omega \times (0, \TM),\\
\tau v_t=\Delta v-v+f(u) &{\rm in}\ \Omega \times (0, \TM),\\
u_\nu=v_\nu=0 &{\rm on}\ \partial\Omega \times (0, \TM),\\
u(x, 0)=u_0(x), \tau v(x,0)= \tau v_0(x) &x \in \bar{\Omega},
\end{cases} 
\end{equation}
and 
\begin{equation}\label{problem_grad}
\tag{$\mathcal{B}$}
\begin{cases}
u_t=\nabla \cdot \left((u+1)^{m_1-1}\nabla u -\chi u(u+1)^{m_2-1}\nabla v\right)+a u^\alpha-b u^\alpha \int_\Omega u^\beta - c |\nabla u|^\delta &{\rm in}\ \Omega \times (0, \TM),\\
\tau v_t=\Delta v-v+f(u) &{\rm in}\ \Omega \times (0, \TM),\\
u_\nu=v_\nu=0 &{\rm on}\ \partial\Omega \times (0, \TM),\\
u(x, 0)=u_0(x), \tau v(x,0)= \tau v_0(x) &x \in \bar{\Omega}.
\end{cases} 
\end{equation}
Herein, $\Omega \subset \mathbb{R}^n$ ($n\in \N$) is a bounded domain with smooth boundary $\partial\Omega$; additionally, $m_1, m_2 \in \R$, $\chi, a, b>0$, $c\geq 0$, $\alpha, \beta, \delta \ge 1$, $\tau \in \{0,1\}$, $u_0(x),  v_0(x)$ are sufficiently regular and nonnegative initial data and $f$ is a suitable smooth function. On the other hand, the subscript $\nu$ in $(\cdot)_\nu$ indicates the outward normal derivative on $\partial\Omega$ and $\TM$ is the maximal existence time up to which solutions to the system are defined.

%From the above discussions, what these models idealize in \emph{chemotaxis} phenomenon is natural; in particular given the abundance of literature on the two models in the absence of logistic terms, or where this is expressed by a polynomial, 
\subsection{Presentation of the main results }
We will address the issue of boundedness of solutions to models \eqref{problem_int} and \eqref{problem_grad}, being our objective imposing conditions on the parameters involved in the gradient-dependent and nonlocal sources, rather than on the coefficients connected to the polynomial part. 

For this reason, it is necessary to specify properties of data describing the problem. In the specific, 
\begin{equation}\label{reglocal}
\begin{cases}
\textrm{ For some } \rho\in(0,1) \textrm{ and } n\in \N, \Omega \subset \R^n \textrm{ is a bounded domain of class } C^{2+\rho}, \textrm{ with boundary } \partial \Omega,\\
   f: [0,\infty) \rightarrow \R^+, \textrm{ with } f \in C^1([0,\infty)), 
   \\ 
    u_0, v_0: \bar{\Omega}  \rightarrow \R^+ , \textrm{ with } u_0,  v_0 \in C_\nu^{2+\rho}(\bar\Omega)=\{\psi \in C^{2+\rho}(\bar{\Omega}): \psi_\nu=0 \textrm{ on }\partial \Omega\}. 
   \end{cases}
\end{equation}
Moreover, for $f$ we might also require that 
\begin{equation}\label{Def_f}
f(s) \leq K s^l \quad \textrm{for } K, l >0.
\end{equation}
%\subsection{Claims of the main results} 
Formally, we will prove the following results, dealing respectively with systems \eqref{problem_int} and \eqref{problem_grad}.
\begin{theorem}\label{MainTheoremInt}
For $\tau\in\{0,1\}$, let $\Omega, f, u_0, v_0$ comply with hypotheses \eqref{reglocal} and \eqref{Def_f}. Additionally, let $m_1, m_2 \in \R$, $\chi, a, b,c>0$, $\beta>\alpha \geq 1$, $\delta \geq \max\{1,m_1\}$. Then, whenever 
\begin{equation*}
 \delta >  \frac{n(m_2-m_1+l)+2(m_2+l)}{2},
\end{equation*} 
%or
%\begin{equation*}
%c>0, \, \beta >m_2+l, \, \delta\geq m_1 \quad \textrm{and} \quad \delta >  \frac{n(\alpha-m_1)+2\alpha }{2},
%\end{equation*} 
problem \eqref{problem_int} admits a unique classical solution $(u,v)=(u(x,t),v(x,t))$, nonnegative and global, in the sense that $u, v \geq 0 $ for all $(x,t) \in \bar{\Omega}\times (0,\infty)$, and uniformly bounded in time, i.e., $u,v \in L^{\infty}(\Omega \times (0,\infty))$.
\end{theorem} 
Some further conditions ensure boundedness also for specific value of $\delta$. We make it precise in the following 
\begin{corollary}\label{corollaryToTheo1}
Under the same hypotheses of Theorem \ref{MainTheoremInt}, there are $k_*>1$, depending on data of model \eqref{problem_int}, and $\mathcal{C}_0=\mathcal{C}_0(k_*), \mathcal{C}_P=\mathcal{C}_P(k_*)>0$ such that the same statement holds 
 if $$\delta = \frac{n(m_2-m_1+l)+2(m_2+l)}{2},$$  whenever
\begin{equation}\label{LimitCaseCInt}
%c > \frac{2^{\delta-1}(k-1)\chi K C_0}{k+m_2-1},\quad \textrm{ if} \quad \tau=0, \quad
c > \frac{2^{\delta-1}(k_*-1)\chi K \mathcal{C}_0}{k_*+m_2-1}\left(\tau \mathcal{C}_P^{\frac{l}{k_*+m_2+l-1}}+(1-\tau)\right).
\end{equation}
\end{corollary}
%%%%%
\begin{theorem}\label{MainTheoremGrad}
For $\tau\in\{0,1\}$, let $\Omega, f, u_0, v_0$ comply with hypotheses \eqref{reglocal} and \eqref{Def_f}. Additionally, let $m_1, m_2 \in \R$, $\chi, a, b>0$, $c\geq0$, $\alpha \geq 1$, $\beta \geq \max\{1,m_1-\alpha\}$ and $1 \leq \delta \leq 2$. Then, whenever either  
\begin{equation*}
c=0, \quad \textrm{and} \quad 
\beta>\max\left\{\frac{n(m_2-m_1+l)+2(m_2+l-\alpha)}{2}, \frac{n(\alpha-m_1)}{2}\right\},
\end{equation*}
or
\begin{equation*}
c>0    \quad \textrm{and} \quad   \beta> \frac{n(\alpha-m_1)}{2} \quad \textrm{and} \quad \delta > \frac{n(m_2+l)}{n+1},
\end{equation*}
or
\begin{equation*}
c>0  \quad \textrm{and} \quad \beta> \frac{n(m_2-m_1+l)+2(m_2+l-\alpha)}{2} \quad \textrm{and} \quad \delta > \frac{n \alpha}{n+1},
\end{equation*}
problem \eqref{problem_grad} admits a unique classical solution $(u,v)=(u(x,t),v(x,t))$, nonnegative and global, in the sense that $u, v \geq 0 $ for all $(x,t) \in \bar{\Omega}\times (0,\infty)$, and uniformly bounded in time, i.e., $u,v \in L^{\infty}(\Omega \times (0,\infty))$.
\end{theorem}
A similar analysis considering more complex situations can be conducted.
\begin{corollary}\label{corollaryToTheo2}
Under the same hypotheses of Theorem \ref{MainTheoremGrad}, and setting $H_{1,\alpha}=1$, for $\alpha=1$, and $H_{1,\alpha}=0$, for $\alpha\neq 1$, there are $M_1>0$, $k_*>1$, depending on data of model \eqref{problem_grad}, and $\mathcal{C}_0=\mathcal{C}_0(k_*), \mathcal{C}_{GN}=\mathcal{C}_{GN}(k_*), \mathcal{C}_P=\mathcal{C}_P(k_*)>0$ such that the same statement holds:
\begin{enumerate}
\item For $c=0$
\begin{itemize}
\item[$\triangleright$] either for 
$\beta=\frac{n(m_2-m_1+l)+2(m_2+l-\alpha)}{2}> \frac{n(\alpha-m_1)}{2}$, provided 
\begin{equation}\label{LimitA}
%b > \frac{2^{\alpha+\beta-2}(k_*-1)\chi K C_0}{k_*+m_2-1}, \quad \textrm{if} \quad \tau=0, \quad 
b > \frac{2^{\alpha+\beta-2}(k_*-1)\chi K \mathcal{C}_0}{k_*+m_2-1} \left(\tau \mathcal{C}_P^{\frac{l}{k_*+m_2+l-1}}+(1-\tau)\right),
\end{equation}
\item[$\triangleright$] or for
$\beta=\frac{n(\alpha-m_1)}{2}> \frac{n(m_2-m_1+l)+2(m_2+l-\alpha)}{2}$ whenever
\begin{equation}\label{LimitB}
a<\frac{b}{\mathcal{C}_0 2^{\alpha +\beta-2}}-\frac{\tau H_{1,\alpha}}{k_*}-b|\Omega|
\quad \textrm{provided} \quad\frac{b}{\mathcal{C}_0 2^{\alpha +\beta-2}}-\frac{\tau H_{1,\alpha}}{k_*}-b|\Omega|>0,
\end{equation}
\item[$\triangleright$] or for $\beta=\frac{n(m_2-m_1+l)+2(m_2+l-\alpha)}{2}= \frac{n(\alpha-m_1)}{2}$ if
\begin{equation}\label{LimitC}
b>2^{\alpha+\beta-2}\mathcal{C}_0\left[\frac{(k_*-1)\chi K}{k_*+m_2-1}\left(\tau \mathcal{C}_P^\frac{l}{k_*+m_2+l-1}+(1-\tau)\right)+a+b|\Omega|+\frac{\tau H_{1,\alpha}}{k_*}\right],
\end{equation}
\end{itemize}
\item For $c>0$
\begin{itemize}
\item[$\triangleright$] either for $\beta>\frac{n(\alpha-m_1)}{2}$ and $\delta=\frac{n(m_2+l)}{n+1}$ if
%\item $c > \frac{(k-1)\chi K}{k+m_2-1}(2 C_{GN})^{\frac{\delta(k+m_2+l-1)}{k-1+\delta}} (M_1+|\Omega|)^{\frac{m_2+l}{n+1}}\left(\frac{k-1+\delta}{\delta}\right)^{\delta}$ if $\tau=0$,
\begin{equation}\label{LimitD}
c > \frac{(k_*-1)\chi K}{k_*+m_2-1}(2 \mathcal{C}_{GN})^{\frac{\delta(k_*+m_2+l-1)}{k_*-1+\delta}} (M_1+|\Omega|)^{\frac{m_2+l}{n+1}}\left(\frac{k_*-1+\delta}{\delta}\right)^{\delta}\left(\tau \mathcal{C}_P^{\frac{l}{k_*+m_2+l-1}}+(1-\tau)\right),
\end{equation}
\item[$\triangleright$] or for $\beta>\frac{n(m_2-m_1+l)+2(m_2+l-\alpha)}{2}$ and $\delta=\frac{n \alpha}{n+1}$ provided 
\begin{equation}\label{LimitD}
c > \left(a+b |\Omega|\right)(2 \mathcal{C}_{GN})^{\frac{\delta(k_*+\alpha-1)}{k_*-1+\delta}} (M_1+|\Omega|)^{\frac{\alpha}{n+1}}\left(\frac{k_*-1+\delta}{\delta}\right)^{\delta}.
\end{equation}
\end{itemize}
\end{enumerate}
\end{corollary}
\begin{remark}[The meaning of the term $\tau H_{1,\alpha}$]\label{RemarkFunzioneH}
It is worthwhile to underline that the term $\tau H_{1,\alpha}$ appearing in \eqref{LimitB} and \eqref{LimitC} vanishes for $\tau=0$ and any $\alpha \geq 1$, or $\tau=1$ and $\alpha>1$. Subsequently, for $\beta=\frac{n}{2}(\alpha-m_1)$, the situation $\tau=\alpha=1$ provides a stronger condition for boundedness. (The reason will be clear in Lemma \ref{LemmaProblemGradEll2}.) Conversely to what said for problem \eqref{problem_grad}, the nature of model \eqref{problem_int}, for which is required that $\beta>\alpha\geq 1$, makes that the value $\alpha=1$ does not have a special role in the analysis: in fact, terms associated to this parameter may be controlled by others involving $\beta$ (see details in Lemma \ref{LemmaProblemIntEll1}).  
\end{remark}
Before proceeding to the proofs of the above claims, we can provide some discussions.
\begin{remark}[Some results on problem \eqref{problem_int} in specific situations]\label{MainTheorem-CorollaryInt}
If we refrain from considering the dissipative impact from the nonlocal term in problem \eqref{problem_int} (essentially for $c=0$), boundendess is achieved for $\beta>m_2+l$ or for $\beta=m_2+l$ provided this largeness requirement on $b$: 
%\begin{equation*}
 %\quad \beta > m_2+l 
%\end{equation*}  
%\begin{itemize}
\begin{equation}\label{RestricionOnbModelNonlocal}
b > \frac{2^{\beta-1}(k_*-1)\chi K}{k_*+m_2-1}\left(\tau \mathcal{C}_P^{\frac{l}{k_*+\beta-1}}+(1-\tau)\right).
\end{equation}
On the other hand, for linear diffusion and sensitivity, as well as linear production, corresponding to take $m_1=m_2=l=1$ in model \eqref{problem_int}, the discussion on the the limit case can be faced in a more explicit way. Indeed, for $\beta=2$ and 
$\alpha=1$, the restriction in \eqref{RestricionOnbModelNonlocal} becomes 
\begin{equation}\label{RestricionOnbModelNonlocalLinear}
b > \frac{(n-2)_+\chi K}{n}\left(\tau \mathcal{C}_P^{\frac{1}{\frac{n}{2}+1}}+(1-\tau)\right).
\end{equation}
This precise expression is connected to the choice $k_*=\frac{n}{2}$; indeed, in the linear scenario whenever $u\in L^\infty((0,\TM);L^k(\Omega))$ for $k>\frac{n}{2}$ one also has that $u\in L^\infty((0,\TM);L^\infty(\Omega))$. (We underline that the above assumption on $b$ is not directly  
obtained from \eqref{RestricionOnbModelNonlocal} with $\beta=2$ for reasons which will be specified later; see $\S$ \ref{RemarkLinear}.) In the specific, for $f(u)=u$ (i.e,. $K=1$),  relation  \eqref{RestricionOnbModelNonlocalLinear} recovers \cite[Theorem 2.5]{TelloWinkParEl} 
and  \cite[Theorem 2.2]{ZhengEtAlLogistiFullyParabolic2018}.
%Consequently, the limit cases studied in the previous corollaries become
%\item (Corollary \ref{MainCorollaryInt3}) $\delta=1+ \frac{l(n+2)}{2}$ 
%yields the constrains
%\begin{equation*}
%c \geq \frac{(k-1)\chi K C_0}{k} \; \text{ if } \tau=0 \quad \text{or} \quad
%c \geq \frac{(k-1)\chi K C_0 C_P^{\frac{l}{k+l}}}{k} \; \text{ if } \tau=1.
%\end{equation*}
\end{remark}

%\begin{remark}
%We observe that the previous results continue to be valid for the linearized versions of problems \eqref{problem_int} and \eqref{problem_grad}, corresponding to take $m_1=m_2=1$.
% in the first equation of such system, we find out that the linear version of Theorem \ref{MainTheoremInt} holds in the case $c=0$ if $\beta>1+l$ or $l<\frac{2}{n}$, and in the case $c>0$ if 
%$\delta>1+ \frac{l(n+2)}{2}$. 
%\end{remark}
\begin{remark}[Some results on problem \eqref{problem_grad} in specific situations]\label{MainTheorem-CorollaryGrad}
As to model \eqref{problem_grad}, we can observe what follows:
\begin{itemize}
\item[$\triangleright$] For $c=0$, problem \eqref{problem_grad} is a generalization of what studied in \cite{chiyoduzgunfrassuviglialoro2024,BianEtAlNonlocal,Latos2020}; in the specific, for $m_1=m_2=l=1$ in Theorem \ref{MainTheoremGrad}, the related results extend \cite[Theorem 2.1]{chiyoduzgunfrassuviglialoro2024} and \cite[Theorem 1]{BianEtAlNonlocal}. Additionally, for $m_1=1$, only, since \cite[Theorem 1, second condition in (9)]{Latos2020} apparently has a gap (see \cite[page 7]{Latos2020}), we herein fix it. 
\item[$\triangleright$] For $c=0$, and $m_1=m_2=1$, 
we have $\beta=\frac{n(\alpha-1)}{2}>\frac{nl+2(1+l-\alpha)}{2}$ and condition 
\eqref{LimitB} is reduced to
\begin{equation*}
a<\frac{b}{\mathcal{C}_0},
\end{equation*}
in which no $|\Omega|$ takes part. (Also in this case such condition is not a direct consequence of \eqref{LimitB}.)
\item[$\triangleright$] For $c>0$, and $m_1=m_2=l=1$, the introduction of the gradient nonlinearity makes that boundedness is achieved for the range of $\alpha$ and $\beta$ larger than those derived in \cite[Theorem 2.1]{chiyoduzgunfrassuviglialoro2024}. 
\end{itemize}
\end{remark}
We conclude this section by comparing the effect of each dissipative source, and by confronting the analysis of the fully parabolic case with that of the parabolic-elliptic one.
\begin{remark}[Nonlocal vs. gradient-dependent damping effects, and the case $\tau=0$ vs. the case $\tau=1$] 
Let us give these comments:
\begin{itemize}
\item[$\triangleright$] A natural question in the context of chemotaxis models with smoothing logistic sources is how their expression effectively provides dissipation to the dynamics. In this sense, in problem \eqref{problem_int} we could compare the dissipative actions connected to $au^\alpha-bu^\beta-c\int_\Omega u^\delta$ and $au^\alpha-bu^\beta-c|\nabla u|^\delta$; apparently, the second source has more equilibrium effects than the first.  Indeed, for the linear case $m_1=m_2=l=1$, as we mentioned above it is known that whenever $\frac{2n}{n+1}<\delta \leq 2$, and 
$\beta>\alpha$,  and independently by the birth rate $\alpha$, for any $c>0$ the source $au^\alpha-bu^\beta-c|\nabla u|^\delta$ suffices to prevent any kind of blow-up. (See \cite[Theorem 1.2]{IshidaLankeitVigliloro-Gradient}.) The associated nonlocal situation (i.e., $\beta>\alpha$, $m_1=m_2=l=1$, and $c>0$ in \eqref{problem_int}), conversely, offers (at least) this less sharp scenario toward boundedness: $\delta>\frac{n+4}{2}>2$.
%, or $\beta>\max\{2,\alpha\}$ and $\delta>\frac{n}{2}(\alpha-1)+\alpha$. 
% as to the other, it is less restrictive only for small values of $\alpha$ (precisely, $\alpha<\frac{n(n+5)}{(n+1)(n+2)}$, i.e., for subquadratic birth rates. 
\item[$\triangleright$] 
As the reader can observe, some constraints appearing in Corollaries \ref{corollaryToTheo1} and \ref{corollaryToTheo2} differ by a power of $\mathcal{C}_P$.  The reason behind this feature will be clear below; essentially, for $\tau=0$ the problems are more manageable because from the second equation one directly has $\Delta v=v-f(u)$. When $\tau=1$, oppositely,  $\Delta v=v-f(u)+v_t$ and parabolic regularity results exactly involving $\mathcal{C}_P$ have to be invoked. 
% Alternatively, one might  %analyze the elliptic case corresponding to $\tau=0$ in a different way, without exploiting the expression of the second equation of both problems.
%More specifically, instead of this, 
%use elliptic regularity results also fot  \textcolor{blue}{CITA} would give similar bounds for our %parameters, in which would appear the constant $C_E$ to the same power of $C_P$.
%\textcolor{blue}{REGOLARITA ELLITTICA, FAI SOLO UN CASO PER DIMOSTRARE CHE COSI I RISULTATI SONO UGUALI}
\end{itemize}
\end{remark}
%\subsubsection{Some results on limit cases for linear models: $m_1,m_2 \in \R$}
%The above theorems continue valid also when limit values of involved parameters are considered, provided some extra restrictions is taken into account. This sequence of corollaries make this more precise.

\section{Some preliminaries and auxiliary tools}\label{PreliminarySection}
Let us now dedicate ourselves to the technical details. We will rely on the following two variants of the Gagliardo--Nirenberg inequality.
\begin{lemma}\label{LemmaGN}
Let $\Omega$ satisfy the hypotheses in \eqref{reglocal}, and let $r\geq 1$  and $0<q\leq p \leq \infty$  fulfilling
$\frac{1}{r} \leq \frac{1}{n} + \frac{1}{p}. $
Then, for 
$\theta= \frac{\frac{1}{q}-\frac{1}{p}}{\frac{1}{q}+\frac{1}{n}-\frac{1}{r}},$
there exists $C_{GN}=C_{GN}(p,q,r,\Omega)>0$ such that for all $\psi \in W^{1,r}(\Omega) \cap L^q(\Omega)$
\[
\|\psi\|_{L^p(\Omega)} \leq C_{GN} (\|\nabla \psi\|_{L^r(\Omega)}^{\theta} \|\psi\|_{L^q(\Omega)}^{1-\theta}+ \|\psi\|_{L^q(\Omega)}).
\]
\begin{proof}
See \cite[Lemma 2.3]{LiLankeitNonLinearity}.
\end{proof}
\end{lemma}
\begin{lemma}\label{GagliardoIneqLemma}
Let $\Omega$ satisfy the hypotheses in \eqref{reglocal}, and let for $n\geq 3$ 
\begin{equation*}%\label{def_of_p}
p:=\frac{2n}{n-2}. 
\end{equation*}
Additionally, let $q, r$ satisfy $1 \le r<q<p$ and $\frac{q}{r}<\frac{2}{r}+1-\frac{2}{p}$. Then for all $\epsilon>0$ there exists $C_0=C_0(\Omega, q, r, \epsilon)>0$ such that for all $\psi \in H^1(\Omega) \cap L^r(\Omega)$, 
\begin{equation}\label{Bian}
\|\psi\|_{L^q(\Omega)}^q \le  C_0 \|\psi\|_{L^r(\Omega)}^\gamma
+\epsilon\|\nabla \psi\|_{L^2(\Omega)}^2+\|\psi\|_{L^2(\Omega)}^2, 
\end{equation}
where 
\begin{equation*}
\lambda:=\frac{\frac{1}{r}-\frac{1}{q}}{\frac{1}{r}-\frac{1}{p}} \in (0,1), \quad \gamma:=\frac{2(1-\lambda)q}{2-\lambda q}.
\end{equation*}
The same conclusion holds for $n\in\{1,2\}$ whenever
$q, r$ fulfill, respectively,  $1 \le r<q$ and $\frac{q}{r}<\frac{2}{r}+2$ and $1 \le r<q$ and $\frac{q}{r}<\frac{2}{r}+1$.
\begin{proof}
The proof is a simple adaptation of \cite[Lemma 3.1]{chiyoduzgunfrassuviglialoro2024}.
\end{proof}
\end{lemma}
We also make use of this result, consequence of Maximal Sobolev Regularity results. 
\begin{lemma}\label{lem:MaxReg}
Let $\Omega$ satisfy the hypotheses in \eqref{reglocal}, and let  $q>1$. Then there is $\mathcal{C}_P=\mathcal{C}_P(\Omega,q)>0$ such that the following holds: Whenever $T\in(0,\infty]$, $I=[0,T)$, $h\in L^q(I;L^q(\Omega))$ and $\psi_0\in W^{2,q}_{\nu}= \{\psi_0\in W^{2,q}(\Omega)\,:\, \partial_\nu \psi_0=0 \,\textrm{ on }\, \partial\Omega\}$, every solution $\psi\in W_{loc}^{1,q}(I;L^q(\Omega))\cap L^q_{loc}(I;W^{2,q}(\Omega))$ of 
 \[
  \psi_t=\Delta \psi- \psi + h\;\; \text{ in }\;\;\Omega\times(0,T);\quad
  \partial_{\nu} \psi=0\;\; \text{ on }\;\;\partial\Omega \times(0,T); \quad \psi(\cdot,0)=\psi_0 \;\; \text{ on }\;\;\Omega  
 \]
 satisfies 
 \[
  \int_0^t  e^s{\int_{\Omega} |\Delta \psi(\cdot,s)|^q}ds \le \mathcal{C}_P \left[\lVert \psi_0 \rVert_{q,1-\frac{1}{q}}^q+\int_0^t e^s {\int_{\Omega} |h(\cdot,s)|^{q}}ds\right] \quad \text{for all } t\in(0,T).
 \]
 \begin{proof}
The proof is based on the classical result in \cite{PrussSchanaubeltMaximalRegul}; for an appropriate adaptation to our case see details, for instance, in \cite{IshidaLankeitVigliloro-Gradient}. 
 \end{proof}
\end{lemma}
\begin{remark}[On the constants $C_{GN}, C_0$ and $\mathcal{C}_P$]\label{RemarkCostanti}
In order to improve the readability and comprehension of our computations we have to make precise what follows:
\begin{itemize} 
\item[$\triangleright$] First, we use this simplification. We will introduce some different $C_{GN}, C_0$ but only one $\mathcal{C}_P$. In particular, since the constants $C_{GN}, C_0, \mathcal{C}_P$ are continuous functions of their arguments, with 
$\mathcal{C}_{GN}, \mathcal{C}_0$ and $\mathcal{C}_P$ in Corollaries \ref{corollaryToTheo1} and \ref{corollaryToTheo2}, we refer to the maximum of those $C_{GN}, C_0$ and $\mathcal{C}_P$ itself.
\item[$\triangleright$]
Secondly, it appears important to specify that the key role of Lemma  \ref{lem:MaxReg} is the existence of the constant $\mathcal{C}_P$, which  remains defined once $n,\Omega$ and $q$ are set. In particular (see \cite[Theorem 2.5]{PrussSchanaubeltMaximalRegul}), $\mathcal{C}_P$ does not depend on the initial configuration $\psi_0$ and the source $h$. As to $\lVert \psi_0 \rVert_{q,1-\frac{1}{q}}$, it represents the norm of $\psi_0$ in the interpolation space $(L^q(\Omega),W^{2,q}_{\nu})_{1-\frac{1}{q},q}.$ (See, for instance, \cite[$\S$1]{LunardiBook}.)  
\end{itemize}
\end{remark}

%\begin{lemma}\label{lem:MaxReg}
% Let $n\in \mathbb{N}$, $\Omega\subset \mathbb{R}^n$ be a bounded and smooth domain and $q\in(1,\infty)$. Then there is $C_P>0$ such that the following holds: Whenever $T\in(0,\infty]$, $I=[0,T)$, $f\in L^q(I;L^q(\Omega))$ and $v_0\in W^{2,q}(\Omega)$ is such that $\partial_{\nu} v_0=0$ on $\partial\Omega$, every solution $v\in W_{loc}^{1,q}(I;L^q(\Omega))\cap L^q_{loc}(I;W^{2,q}(\Omega))$ of
%\[
% v_t=\Delta v-v + f\;\; \text{ in }\;\;\Omega\times(0,T);\quad
% \partial_{\nu} v=0\;\; \text{ on }\;\;\partial\Omega \times(0,T); \quad v(\cdot,0)=v_0 \;\; \text{ on }\;\;\Omega
%\]
%satisfies
%\[
% \int_0^t e^s \left(\int_\Omega |\Delta v(\cdot,s)|^q\right)ds \le C_P \left[1+\int_0^t e^s \left(\int_\Omega |f(\cdot,s)|^{q}\right)ds\right] \quad \text{for } 0<t<T.
%\]
%\end{lemma}
This ODE comparison principle will be as well employed.
\begin{lemma}\label{LemmaODI-Comparison}
Let\/ $T>0$ and $\phi:(0,T)\times \R^+_0\rightarrow \R$. If $0\leq y\in C^0([0,T))\cap  C^1((0,T))$ is such that 
\begin{equation*}
y'\leq \phi(t,y)\quad \textrm{for all } t \in (0,T), 
\end{equation*}
and there is $y_1>0$ with the property that whenever $y>y_1$ for some $t\in (0,T)$ one has that $\phi(t,y)\leq 0$, then
\begin{equation*}
y\leq \max\{y_1,y(0)\}\quad \textrm{on } (0,T).
\end{equation*}
\begin{proof}
See \cite[Lemma 3.3]{chiyoduzgunfrassuviglialoro2024}.
\end{proof}
\end{lemma}
Further, we need to adjust parameters taking part in our derivations.
\begin{lemma}\label{ThetaSigma}
Let $n \in \N$, $m_1, m_2 \in \R$, $l >0$ and $\alpha, \delta \geq 1$. Then, there exists $k_*>1$ such that, for all $k>k_*$ and
\begin{equation*}
\theta_1(k)\coloneqq\frac{\frac{k+m_1-1}{2}-\frac{k+m_1-1}{2k}}{\frac{k+m_1-1}{2}-\frac{1}{2}+\frac{1}{n}}, \quad 
\sigma_1(k)\coloneqq \frac{2k}{k+m_1-1}, \quad 
\theta_2(k)\coloneqq\frac{\frac{k-1+\delta}{\delta}-\frac{k-1+\delta}{\delta(k+\alpha-1)}}{\frac{k-1+\delta}{\delta}-\frac{1}{\delta}+\frac{1}{n}}, 
\end{equation*}
\begin{equation*}
\sigma_2(k)\coloneqq\frac{\delta(k+\alpha-1)}{k-1+\delta}, \quad
\theta_3(k)\coloneqq\frac{\frac{k-1+\delta}{\delta}-\frac{k-1+\delta}{\delta(k+m_2+l-1)}}{\frac{k-1+\delta}{\delta}-\frac{1}{\delta}+\frac{1}{n}}, \quad
\sigma_3(k)\coloneqq\frac{\delta(k+m_2+l-1)}{k-1+\delta},
\end{equation*}
\begin{equation*} 
\theta_4(k)\coloneqq\frac{\frac{k+m_1-1}{2}-\frac{1}{2}}{\frac{k+m_1-1}{2}-\frac{1}{2}+\frac{1}{n}}, 
\end{equation*}
these relations hold:
\begin{table}[H]
\centering
\begin{subequations}
\begin{subtable}[h]{0.24\textwidth}
\centering
\begin{equation}\label{theta1}
0<\theta_1<1,
\end{equation}	
\end{subtable}
\hfill
\begin{subtable}[h]{0.24\textwidth}
\centering
\begin{equation}\label{theta2}
0<\theta_2<1,
\end{equation}	
\end{subtable}
\begin{subtable}[h]{0.24\textwidth}
\centering
\begin{equation}\label{theta3}
0<\theta_3<1,
\end{equation}	
\end{subtable}
\hfill
\begin{subtable}[h]{0.24\textwidth}
\centering
\begin{equation}\label{theta4}
0<\theta_4<1,
\end{equation}	
\end{subtable}
\end{subequations}
\end{table}
\begin{subequations}
\begin{equation}\label{sigmatheta1}
\frac{\sigma_1\theta_1}{2}>0,
\end{equation}
\begin{equation}\label{sigmatheta2}
0<\frac{\sigma_2\theta_2}{\delta}<1 \quad \textrm{if} \quad \delta > \frac{n \alpha}{n+1},
\end{equation}	
\begin{equation}\label{sigmatheta3}
0<\frac{\sigma_3\theta_3}{\delta}<1 \quad \textrm{if} \quad \delta > \frac{n(m_2+l)}{n+1}.
\end{equation}	
\end{subequations}
\begin{proof}
The expressions defining $\theta_1(k), \theta_2(k), \theta_3(k)$ and $\theta_4(k)$ can equivalently be expressed as $\frac{a_1(k) - \frac{a_1(k)}{a_2(k)}}{a_1(k) - a_3}$, with certain $a_1(k),a_2(k) \in O(k)$ and $a_3\in\mathbb{R}$. Let us choose $\theta_1(k)$ as example while proving \eqref{theta1}. For \eqref{theta2}, \eqref{theta3} and \eqref{theta4} the reasonings are similar. First of all, we rewrite $\theta_1(k)$ as
\begin{equation*}
\theta_1(k) = \frac{\frac{k+m_1 - 1}{2} - \frac{k+m_1 -1}{2k}}{\frac{k+m_1 - 1}{2} + \frac1n - \frac12} = \frac{ 1 - \frac{1}{k}}{1 - \frac{n-2}{n(k+m_1-1)}}.
\end{equation*}
Therefore, $\lim_{k\to\infty} \theta_1(k) = 1$ and, as a consequence, $\theta_1 > 0$ definitively. To prove that $\theta_1 <1$, we notice that
\begin{equation*}
\frac{1}{k} > \frac{n-2}{n(k+m_1-1)} \Leftrightarrow \frac{k+m_1-1}{k} > \frac{n-2}{n} = 1 - \frac2n,
\end{equation*} 
and this is a natural implication of what follows:
\begin{equation*}
 \lim_{k\to\infty} \frac{k+m_1-1}{k} = 1 \Rightarrow  \frac{k+m_1-1}{k} > 1 - \frac2n, \quad \textrm{definitively}.
\end{equation*}
Therefore, $\theta_1 \in (0,1)$ for $k$ sufficiently large.

As to \eqref{sigmatheta1}, \eqref{sigmatheta2} and \eqref{sigmatheta3}, since $\sigma_1(k), \sigma_2(k)$, and $\sigma_3(k)$ are definitively positive, the same applies for $\sigma_i \theta_i$, with $i=1,2,3$. Let us prove only \eqref{sigmatheta2}: we notice that
\begin{equation*}
\frac{\sigma_2(k)\theta_2(k)}{\delta} = 1 - \frac{\delta(1+\frac1n)-\alpha}{k-2+\delta(1+\frac{1}{n})}\in (0,1) \Leftrightarrow 
\frac{\delta(1+\frac1n)-\alpha}{k-2+\delta(1+\frac{1}{n})}\in (0,1).
\end{equation*}
Being $\lim_{k\to\infty} \frac{\delta(1+\frac1n)-\alpha}{k-2+\delta(1+\frac{1}{n})}= 0$, the inclusion holds for $\delta>\frac{n\alpha}{n+1}$ as claimed. 

As a consequence, a large $k_*>1$ such that our statements are satisfied for all $k>k_*$ can be detected. 
\end{proof}
\end{lemma}
The following consideration is essential for understanding most of the remaining reasonings.  
\begin{convention}[On the role of $k_*$ in Corollaries \ref{corollaryToTheo1} and \ref{corollaryToTheo2}]\label{Convenzione} For reasons which will be clear later on, we are interested in detecting the minimum real greater than $1$, in a such a way that all the derivations in this paper remain justified for larger values of this minimum. In particular, in some forthcoming analyses we might require to progressively increase the value $k_*>1$ introduced in Lemma \ref{ThetaSigma}. We will tacitly carry out this procedure and, without any relabeling, with $k_*$ we might refer to a higher value. 
\end{convention}
\section{Local solutions and their main properties: a boundedness criterion}\label{SectionLocalInTime}
\begin{lemma}[Local existence]\label{localSol}
For $\tau\in \{0,1\}$, let $\Omega$, $f, u_0, v_0$ comply with hypotheses in \eqref{reglocal}.  Additionally, let $m_1, m_2 \in \R$, $\chi, a, b>0$, $c\geq 0$ and $\alpha, \beta, \delta \geq 1$.  
Then problems \eqref{problem_int} and \eqref{problem_grad} have a unique and nonnegative classical solution
\begin{equation*}
(u,v)\in C^{2+\rho,1+\frac{\rho}{2}}( \Bar{\Omega} \times [0, \TM))\times C^{2+\rho,\tau+\frac{\rho}{2}}( \Bar{\Omega} \times [0, \TM)),
\end{equation*}
for some maximal $\TM\in(0,\infty]$, obeying the following property:
 \begin{equation}\label{dictomyCriteC2+delVerions0}
 \text{either } \TM=\infty \quad \text{or}\quad \limsup_{t \to \TM} \left(\|u(\cdot,t)\|_{C^{2+\rho}(\bar\Omega)}+\|v(\cdot,t)\|_{C^{2+\rho}(\bar\Omega)}\right)=\infty.
\end{equation}
% If, additionally $\delta\leq2$ in problem \eqref{problem_grad}, then for both systems $\TM$ is such that 
%\begin{equation}\label{dictomyCriteC2+del}
% \text{either } \TM=\infty \quad \text{or}\quad \limsup_{t\to \TM} \lVert u(\cdot,t)\rVert_{L^{\infty}(\Omega)} = \infty.
%\end{equation}
Finally, if $\beta>\alpha$
\begin{equation}\label{MassBounded_int}
\int_\Omega u(x,t)\,dx \leq M_0:= \max\left\{\int_\Omega u_0(x)dx, \left(\frac{a}{b}|\Omega|^{\beta-\alpha}\right)^\frac{1}{\beta-\alpha}\right\}
\quad \textrm{for all } t \in (0,\TM),
\end{equation}
for problem \eqref{problem_int} and 
\begin{equation}\label{MassBounded_grad}
\int_\Omega u(x,t)\,dx \leq M_1:=\max\left\{\int_\Omega u_0(x)dx, \left(\frac{a}{b |\Omega|^{1-\beta}}\right)^{\frac{1}{\beta}}\right\} 
\quad \textrm{on } (0,\TM),
\end{equation}
for problem \eqref{problem_grad}.
\begin{proof} 
Letting $X=(x,t)$ and %$B(X,u,\nabla u)=F(X,u)-c|\nabla u|^\gamma$ with
\begin{equation}\label{ExpressionofBForSources}
B(X,u,\nabla u)=B=
\begin{cases}
au^\alpha-b u^\beta-c\int_\Omega u^\delta, \\
au^\alpha-b u^\alpha \int_\Omega u^\beta-c|\nabla u|^\delta, 
\end{cases}
\end{equation}
problems \eqref{problem_int} and \eqref{problem_grad} can essentially be collected in this single system:
\begin{equation*}%\label{problem_Single}
\begin{cases}
u_t=\nabla \cdot \left((u+1)^{m_1-1}\nabla u -\chi u(u+1)^{m_2-1}\nabla v\right)+B(X,u,\nabla u)&{\rm in}\ \Omega \times (0, \TM),\\
\tau v_t=\Delta v-v+f(u) &{\rm in}\ \Omega \times (0, \TM),\\
u_\nu=v_\nu=0 &{\rm on}\ \partial\Omega \times (0, \TM),\\
u(x, 0)=u_0(x), \tau v(x,0)= \tau v_0(x) &x \in \bar{\Omega}.
\end{cases} 
\end{equation*}
Let us give some hints (more details are available in \cite[Lemma 2.1]{IshidaLankeitVigliloro-Gradient} for $B$ depending on $\nabla u$ and in \cite[Proposition 4]{BianEtAlNonlocal} for the other situation) only concerning the existence issue.  
%and $0\not\equiv u_0, \tau v_0,\tau w_0\in C_\nu^{2+\delta}(\bar\Omega)$ and $0\not\equiv v_0\in W^{1,%\infty}(\Omega)$, such that $\lVert u_0 \rVert_{L^\infty(\Omega)}\leq R$, 

For any $R>0$, let us consider for some $0<T\leq 1$, to be determined, the closed convex
subset
$$S_{T}=\{0\le u \in C^{1,\frac{\rho}{2}}(\bar{\Omega}\times [0,T]):  \|u(\cdot,t)-u_0\|_{L^\infty(\Omega)}\le R,\, \; \text{for\ all}\
t\in[0,T]\}.$$
Once an element $\tilde{u}$ of $S_T$ is picked, from the properties of $f$ one has that $f(\tilde{u})\in C^{\rho,\frac{\rho}{2}}(\bar{\Omega}\times [0,T])$, so that the related solution $v$ to problem
 \begin{equation}\label{2.2}
		\begin{cases}
	\tau v_t-	\Delta v+ v=f(\tilde{u})&\text{in}\ \Omega\times(0,T),\\
		v_\nu=0&\text{on}\ \partial\Omega\times(0,T),\\
		\tau v(x,0)=\tau v_0(x),
	\end{cases}
\end{equation}
%and
%\begin{equation}\label{2.3}
%	\begin{cases}
%	\tau w_t-	\Delta w+ w=f_2(\tilde{u})&\text{in}\ \Omega\times(0,T),\\
%		w_\nu=0&\text{on}\ \partial\Omega\times(0,T),\\
%		\tau w(x,0)=\tau w_0(x),
%	\end{cases}
%\end{equation}
is, thanks to $\tau v_0 \in C_\nu^{2+\rho}(\bar{\Omega})$ and  classical elliptic and parabolic regularity results (\cite[Theorem 9.33]{BrezisBook},  \cite[Theorem 6.31]{GilbarTrudinger} for $\tau=0$, and \cite[Corollary 5.1.22]{LunardiBook}, \cite[Theorem IV.5.3]{LSUBookInequality} for $\tau=1$), such that 
\begin{equation*}
v\in C^{2+\rho, \tau+\frac{\rho}{2}}(\bar\Omega \times [0,T])
\quad \textrm{and precisely}\quad \sup_{t\in [0,T]} \lVert v(\cdot, t)\rVert_{C^{2+\rho}(\bar\Omega)}\leq H,
\end{equation*}
where $H=H(f,\tilde{u},\tau v_0)>0$ is estimated by $\tilde{u}$, and henceforth $H=H(R)$.
%subsequently,  one gets that $\nabla v,\nabla w\in L^\infty(\Omega)$ for all $t\in(0,T)$, being the norm .
Now, due to obtained properties of $\nabla v$ and the smoothness of $\zeta \longmapsto (1+\zeta)^\vartheta$ for all $\vartheta \in \R$ and $\zeta\geq 0$, let us consider \begin{equation}\label{2.1}
	\begin{cases}
		u_t=\nabla\cdot\left((\tilde{u}+1)^{m_1-1}\nabla u-\chi u(\tilde{u}+1)^{m_2-1}\nabla v
		\right) +B(X,u,\nabla u)
		&\text{in}\ \Omega\times(0,T),\\
		u_\nu=0&\text{on}\ \partial\Omega\times(0,T),\\
		u(x,0)=u_0(x)&x\in\overline{\Omega}.
	\end{cases}
\end{equation}
In the same spirit of \cite{IshidaLankeitVigliloro-Gradient} and \cite{BianEtAlNonlocal}
%\begin{equation}\label{PositionForLunardi}
%\begin{cases}
%\varphi_1: \Omega \times (0,T)\rightarrow \R^+, \textrm{ with }\varphi_1(x,t):= (\tilde{u}(x,t)+1)^{m_1-1},\\
%\varphi_2: \Omega \times (0,T)\rightarrow \R^n, \textrm{ with }\varphi_2(x,t):= -\chi (\tilde{u}(x,t)+1)^{m_2-1}\nabla v(x,t)+  \xi (\tilde{u}(x,t)+1)^{m_3-1}\nabla w(x,t),\\
%\varphi_3: \Omega \times (0,T)\times \R^{1+n}\rightarrow \R, \textrm{ with }\varphi(x,t,u, \nabla u):= \lambda u^\rho(x,t)-\mu u^k(x,t)-c|\nabla u(x,t)|^\gamma,
%\end{cases}
%\end{equation}
%and also this system: 
%In the specific, with the positions in \eqref{PositionForLunardi}, problem \eqref{2.1} reads as
%\begin{equation*}%\label{2.1BIS}
%	\begin{cases}
%		u_t=\mathcal{A}u +\varphi(x,t,u,\nabla u)=\varphi_1 \Delta u +u \nabla \cdot \varphi_2+(\nabla %\varphi_1+\varphi_2)\cdot \nabla u + \varphi_3
%		&\text{in}\ \Omega\times(0,T),\\
%		\mathcal{B}_1u=u_\nu=0&\text{on}\ \partial\Omega\times(0,T),\\
%		u(x,0)=u_0(x)&x\in\overline{\Omega},
%	\end{cases}
%\end{equation*}
%and it is seen that $\varphi$ comply with \cite[(7.3.5)]{LunardiBook}, exactly using the notation therein. 
by recalling that $u_0\in C^{2+\rho}_\nu(\bar\Omega)$, it can be seen that there exists some $0<T< 1$, such that problem \eqref{2.1} has a unique solution 
\begin{equation*}
u\in C^{2+\rho,1+\frac{\rho}{2}}(\bar\Omega\times [0,T]).
\end{equation*} 
In particular,  this produces some positive constant $L=L(\tilde{u},\nabla v,u_0)$, i.e., $L=L(R)$, with the property that
$$|u(x,t)-u_0(x)|\le Lt^{1+\tfrac{\rho}{2}} \quad \textrm{for all} \quad(x,t)\in\Omega\times(0,T), \quad \textrm{or}\quad \max_{t\in[0,T]}\,\|u(\cdot,t)-u_0\|_{L^\infty(\Omega)}\le LT^{1+\frac{\rho}{2}}. $$
In this way, we deduce %for $T\leq \left(\frac{R}{K}\right)^\frac{2}{2+\delta}$,
that 
$$ \textrm{for} \quad T\leq \left(\frac{R}{L}\right)^\frac{2}{2+\rho}, \quad\|u(\cdot,t)-u_0\|_{L^\infty(\Omega)}\le R\;\quad
\textrm{for all} \quad t\in[0,T]. $$
Moreover, since $\underline{u}\equiv 0$ is a subsolution of the first equation in \eqref{2.1}, the parabolic comparison principle warrants the nonnegativity of $u$ on $\Omega \times (0,T)$. So, for values of $T$ as above, the map  $\Phi (\tilde{u})=u$ where $u$ solves problem \eqref{2.1}, is such that $\Phi(S_T) \subset S_T$ and $\Phi$ is compact, because the
 \cite[Ascoli--Arzel\`a Theorem 4.25]{BrezisBook} implies that the natural embedding
of $C^{2+\rho,\tau+\frac{\rho}{2}}(\bar{\Omega}\times [0,T]))$ into $C^{1,\frac{\rho}{2}}(\bar{\Omega}\times [0,T]))$ is a compact linear operator. Let $u$ be the fixed point of $\Phi$ asserted by the Schauder's fixed point theorem: first, the elliptic and parabolic maximum principles and $f(u)\geq 0$ in problems \eqref{2.2} also imply $v\geq 0$ in $\Omega \times (0,T)$; secondly,  as seen, $u,v$ have the claimed regularity. 

On the other hand, taking $T$ as initial time and $u(\cdot, T)$ as the initial condition, the above explained procedure would provide a solution $(\hat{u},\hat{v},\hat{w})$ defined on $\bar{\Omega}\times [T,\hat{T}]$, for some $\hat{T}>0$, which by uniqueness would be the prolongation of $(u,v,w)$ (exactly, from $\bar{\Omega}\times [0,T]$ to  $\bar{\Omega}\times [0,\hat{T}])$. This machinery can be repeated up to construct a maximal interval time $[0,\TM)$ of existence, in the sense that either $\TM=\infty$, or if $\TM<\infty$ no solution belonging to $C^{2+\rho,\tau+\frac{\rho}{2}}(\bar{\Omega}\times [0,\TM])$ may exist and henceforth relation \eqref{dictomyCriteC2+delVerions0} has to be fulfilled.

On the other hand, as to relation \eqref{MassBounded_int}, an integration of the first equation of problem \eqref{problem_int} yields
\begin{equation}\label{M01}
\begin{split}
y':=\frac{d}{dt} \int_\Omega u= a \int_\Omega u^\alpha - b \int_\Omega u^\beta - c \int_\Omega u^\delta \leq a \int_\Omega u^\alpha - b \int_\Omega u^\beta \quad \textrm{for all } t \in (0,\TM).
\end{split}
\end{equation}
%Since $\beta>\alpha$, an application of the H\"{o}lder and Young's inequalities to the first integral of the right-hand side of \eqref{M01} leads to
%\begin{equation*} 
% a \int_\Omega u^\alpha  \leq  \left(\int_\Omega u^\beta \right)^{\frac{\alpha}{\beta}} 
% |\Omega|^{\frac{\beta-\alpha}{\beta}} \leq \frac{b}{2} \int_\Omega u^\beta + \frac{\beta-\alpha}{\beta} \left(\frac{b\beta}{2a\alpha}\right)^{-\frac{\alpha}{\beta-\alpha}}|\Omega| \quad \textrm{on } (0,\TM),
%\end{equation*}
%which inserted in the previous estimate gives
%\begin{equation}\label{M0S}
%\frac{d}{dt} \int_\Omega u \leq -\frac{b}{2} \int_\Omega u^\beta + \frac{\beta-\alpha}{\beta} \left(\frac{b\beta}{2a\alpha}\right)^{-\frac{\alpha}{\beta-\alpha}}|\Omega| \quad \textrm{for all } t \in (0,\TM).
%\end{equation}
Now, thanks to the H\"{o}lder inequality, we have for $\beta>\alpha\geq 1$
\begin{equation}\label{uk}
    -\into u^\beta \leq -\abs*{\Omega}^{\frac{\alpha-\beta}{\alpha}}\tonda*{\into u^\alpha}^\frac{\beta}{\alpha} \quad \textrm{and} \quad -\into u^\alpha \leq -\abs*{\Omega}^{1-\alpha}\tonda*{\into u}^\alpha\quad \text{in }(0,\TM).
\end{equation}
By combining expressions \eqref{M01} and \eqref{uk} we arrive for $\gamma(t)=a \int_\Omega u^\alpha(x,t)dx\geq 0$ on $(0,\TM)$ at this initial problem
\begin{equation*}
\begin{dcases}
    y'\leq \gamma(t)\left(1 - \frac{b}{a} \abs*{\Omega}^{\alpha-\beta} y^{\beta-\alpha}\right) & \text{in }(0,\TM),\\
    y(0)=\into u_0(x)dx,
\end{dcases}
\end{equation*}
so concluding by invoking Lemma \ref{LemmaODI-Comparison} with $$T=\TM,\, \phi(t,y)=\gamma(t)\left(1 - \frac{b}{a} \abs*{\Omega}^{\alpha-\beta} y^{\beta-\alpha}\right),\, y(0)=\into u_0(x)dx\, \textrm{ and }\, y_1=\left(\frac{a}{b} |\Omega|^{\beta-\alpha}\right)^{\frac{1}{\beta-\alpha}}.$$
%From a further application of the H\"{o}lder inequality, we deduce that
%\begin{equation*} 
%\int_\Omega u \leq \left(\int_\Omega u^\beta \right)^{\frac{1}{\beta}} 
% |\Omega|^{\frac{\beta-1}{\beta}} \quad \textrm{for all } t \in (0,\TM).
%\end{equation*}
%After some manipulations in the previous bound, by inserting it in \eqref{M0S}, we obtain
%\begin{equation*}
%\frac{d}{dt} \int_\Omega u \leq -\frac{b}{2}  |\Omega|^{1-\beta} \left(\int_\Omega u\right)^{\beta} +  \frac{\beta-\alpha}{\beta} \left(\frac{b\beta}{2a\alpha}\right)^{-\frac{\alpha}{\beta-\alpha}}|\Omega| \quad \textrm{on } (0,\TM),
%\end{equation*}
%henceforth, by ODI comparison principles it follows the claim.
For \eqref{MassBounded_grad}, similarly to what done, by integrating the first equation of problem \eqref{problem_grad} we have 
\begin{equation*}
\frac{d}{dt} \int_\Omega u \leq \int_\Omega u^\alpha \left(a- b \int_\Omega u^\beta \right) \quad \textrm{for all } t \in (0,\TM).
\end{equation*}
At this point, by reasoning as before (details available in \cite[Lemma 4.1]{chiyoduzgunfrassuviglialoro2024}), we conclude by invoking again Lemma \ref{LemmaODI-Comparison}.
\end{proof}
\end{lemma}
From now on with $(u,v)$ we will refer to the unique local solutions to models \eqref{problem_int} and \eqref{problem_grad} defined on $\Omega \times (0,\TM)$ and  provided by Lemma \ref{localSol}. 
As seen in such a Lemma, relation \eqref{dictomyCriteC2+delVerions0} holds and in this case the solution is said to blow-up in finite time $\TM$ in $C^{2+\rho}(\Omega)$-norm, and also in $L^\infty(\Omega)$-norm for problem \eqref{problem_int}. Conversely, in order to have this implication for model \eqref{problem_grad}, where gradient-dependent sources appear, 
%exactly in order to rephrase the extensibility criterion in \eqref{dictomyCriteC2+del} in terms of some uniform-in-time boundendess of $u$, 
we have to assume some restriction on the growth of $|\nabla u|^\delta$ for $\delta \geq 1$.
% to those whose growth is at the most quadratic. 
We precisely have this 
\begin{lemma}[Boundedness criterion] \label{ExtensionLemma}
If for $\delta\geq 1$ ($1\leq \delta\leq 2$) the solution $(u,v)$ to problem \eqref{problem_int} (problem \eqref{problem_grad}) is such that $u\in L^\infty((0,\TM);L^\infty(\Omega))$, then we have $\TM=\infty$ and in particular $u\in L^\infty((0,\infty);L^\infty(\Omega))$. 
\begin{proof}
The claim is achieved by contradiction. Let $\TM$ be finite; it is sufficient to show that the uniform-in-time boundedness of $u$ in $L^\infty(\Omega)$ for $t\in(0,\TM)$, entails uniform-in-time boundedness of $\lVert u(\cdot,t) \rVert_{C^{2+\rho}(\bar{\Omega})}+\lVert v(\cdot,t) \rVert_{C^{2+\rho}(\bar{\Omega})}$ on $[0,\TM].$ In fact, from \eqref{dictomyCriteC2+delVerions0}, we would have an inconsistency. 
%\begin{equation}
%\text{if} \quad \TM<\infty, \quad \text{then} \quad \lim_{t \to \TM} \left(\|u(\cdot,t)\|_{C^{2+\delta}(\bar\Omega)}+\|v(\cdot,t)\|_{C^{2+\delta}(\bar\Omega)}+\|w(\cdot,t)\|_{C^{2+\delta}(\bar\Omega)}\right)=\infty.
%\end{equation}

On the one hand, we have boundedness of $\lVert v(\cdot,t) \rVert_{C^{1+\rho}(\bar{\Omega})}$ on $[0,\TM]$, as consequence of our hypothesis $u\in L^\infty((0,\TM);L^\infty(\Omega))$, once \cite[Theorem~1.2]{lieberman_paper} (if $\tau=1$) or \cite[Theorem.~I.19.1]{friedman_pde} (if $\tau=0$) and the Sobolev embedding immersion $W^{2,k}(\Omega)\subset C^{1+\rho}(\bar{\Omega})$, with arbitrarily large $k$, are invoked. 

Now for $X=(x,t)$, $A(X,u,\nabla u)=(u+1)^{m_1-1}\nabla u-\chi u(u+1)^{m_2-1}\nabla v(x,t) \textrm{ and } B(X,u,\nabla u)$ as in \eqref{ExpressionofBForSources}
we consider this system:
\begin{equation}\label{ProblemAuxiliary}
\begin{cases}
u_t=\nabla \cdot A(X,u,\nabla u)+B(X,u,\nabla u) & \textrm{ in } \Omega \times (0,\TM),\\
A(X,u,\nabla u)\cdot \nu=0 & \textrm{ on } \partial \Omega \times (0,\TM)\\
u(x,0)=u_0(x)& x\in \bar{\Omega}.
\end{cases}
\end{equation}
For $B$ not depending on $\nabla u$, by using that $u$ is a uniform-in-time bounded solution to problem \eqref{ProblemAuxiliary},  we indeed have from parabolic regularity results that $\lVert u(\cdot,t) \rVert_{C^{1+\rho}(\bar{\Omega})}$ is finite as well on $[0,\TM]$. From this higher regularity of $u$,  bootstrap  arguments already developed in Lemma \ref{localSol}, in conjunction with the regularity of the initial data $(u_0,\tau v_0)$ for model \eqref{problem_int} give for any $\delta \geq 1$
$$
\sup_{t\in [0,\TM]}\left(\lVert u(\cdot,t) \rVert_{C^{2+\rho}(\bar{\Omega})}+\lVert v(\cdot,t) \rVert_{C^{2+\rho}(\bar{\Omega})}\right)<\infty.
$$ 
As to problem \eqref{problem_grad}, for $\delta\in [1,2]$ the gradient-dependent term $B$ behaves as in 
\cite[($1.2c$),  Theorem.~1.2]{lieberman_paper}; subsequently, we invoke \cite[Theorem.~1.2]{lieberman_paper} to see that $\lVert u(\cdot,t) \rVert_{C^{1+\rho}(\bar{\Omega})}$ is finite as well on $[0,\TM]$, so similarly concluding. 
\end{proof}
\end{lemma} 
\section{Achieving uniform-in-time boundedness of $u$ in appropriate Lebesgue spaces. Proof of the results}\label{EstimatesAndProofSection}
Since, as we will discuss, the uniform-in-time boundedness of $u$ is implied whenever $u\in L^\infty((0,\TM);L^k(\Omega))$ for some $k>k_*$, here under we derive some \textit{a priori} integral estimates. 

Let us define these functionals 
\begin{equation}\label{Func}
y(t):=\int_\Omega (u+1)^k \quad \textrm{for all } t \in (0,\TM),
\end{equation}
and 
\begin{equation}\label{Def_Fm2}
F_{m_2}(u):= \int_0^u \hat{u}(\hat{u}+1)^{k+m_2-3}\,d\hat{u} \quad \textrm{on } (0,\TM),
\end{equation}
this one naturally fulfilling   
\begin{equation}\label{Fm2}
0\leq F_{m_2}(u)\leq \frac{1}{k+m_2-1}[(u+1)^{k+m_2-1}-1] \quad \textrm{for all } t \in (0,\TM).
\end{equation}
The next two results are inequalities valid for a class of certain functions, and obtained by exploiting Lemma \ref{LemmaGN} and Lemma \ref{GagliardoIneqLemma}. Nevertheless, we will directly derive such inequalities for the component $u$ of the solution to our models. 

First we point out that
\begin{itemize}
\item [$\triangleright$] all the constants $c_i$, $i = 1, 2, \ldots$, appearing below, are assumed to be positive;
\item [$\triangleright$] with $\epsilon$ we indicate an arbitrary positive real number, and the multiplication by another constant and the sum with another arbitrary constant are not performed, and the final result is for commodity as well labeled with $\epsilon$;
\item [$\triangleright$] $C_{GN}, C_0, \mathcal{C}_P$ and $M_0, M_1$ are the constants defined in Lemmas \ref{LemmaGN}, \ref{GagliardoIneqLemma}, \ref{lem:MaxReg} and \ref{localSol}, respectively. 
\end{itemize} 
Let us start taking advantage from Lemma \ref{LemmaGN}.
\begin{lemma}\label{GNandInterp}
Let the assumptions of Lemma \ref{ThetaSigma} be valid and let $k_*>1$ be the value described in Convention \ref{Convenzione}. Then for every $\hat{c}>0$, $k>k_*$ and all $t\in (0,\TM)$ we have
\begin{equation}\label{GN2}
-\int_\Omega |\nabla (u+1)^{\frac{k+m_1-1}{2}}|^2 \leq \const{gn21} - \const{gn22} 
\left(\int_\Omega (u+1)^k\right)^{\frac{k+m_1-1}{k \theta_1}}  \quad 
\end{equation}

\begin{equation}\label{GN3}
\hat{c} \int_\Omega (u+1)^{k+\alpha-1} \leq
\begin{cases}
\displaystyle  \frac{kc}{2} \left(\frac{k-1+\delta}{\delta}\right)^{-\delta} \int_\Omega |\nabla (u+1)^{\frac{k-1+\delta}{\delta}}|^\delta+\const{gn3} & \textrm{ if } \displaystyle   \delta > \frac{n \alpha}{n+1},\\
\displaystyle    \hat{c}\, (2 C_{GN})^{\sigma_2} (M_1+|\Omega|)^{\frac{\alpha}{n+1}} \int_\Omega |\nabla (u+1)^{\frac{k-1+\delta}{\delta}}|^\delta+\const{gn03} & \textrm{ if } \displaystyle   \delta = \frac{n \alpha}{n+1},
 \end{cases}
\end{equation}

\begin{equation}\label{GN4}
\begin{split}
\hat{c} \int_\Omega (u+1)^{k+m_2+l-1}
&\leq 
\begin{cases}
\displaystyle \frac{kc}{2} \left(\frac{k-1+\delta}{\delta}\right)^{-\delta} \int_\Omega |\nabla (u+1)^{\frac{k-1+\delta}{\delta}}|^\delta+\const{gn4}
 & \textrm{ if } \displaystyle  \delta >\frac{n(m_2+l)}{n+1},\\
\displaystyle  \hat{c}\, (2 C_{GN})^{\sigma_3} (M_1+|\Omega|)^{\frac{m_2+l}{n+1}} \int_\Omega |\nabla (u+1)^{\frac{k-1+\delta}{\delta}}|^\delta+\const{gn04}
 & \textrm{ if } \displaystyle  \delta =\frac{n(m_2+l)}{n+1}.
 \end{cases}
 \end{split}
\end{equation}
\begin{proof}
In the sequel we will use (possibly without mentioning it if not necessary)
the well-known inequality 
\begin{equation}\label{ABIn}
(A+B)^s \leq c(s) (A^s+B^s) \quad \text{for all } A,B \geq 0 \; \text{and }
c(s):=
\begin{cases}
2^s &\text{if } s> 0,\\
2^{s-1} &\text{if } s\geq 1.
\end{cases}
\end{equation}
Let us show \eqref{GN2}. By recalling \eqref{theta1} and \eqref{sigmatheta1}, the boundedness of the mass \eqref{MassBounded_int} for problem \eqref{problem_int} and 
 \eqref{MassBounded_grad} for problem \eqref{problem_grad}, and the Gagliardo--Nirenberg and the Young inequalities, provide for all $t \in (0,\TM)$
\begin{equation*}
\begin{split}
\int_\Omega (u+1)^k&=\|(u+1)^{\frac{k+m_1-1}{2}}\|_{L^{\frac{2k}{k+m_1-1}}(\Omega)}^{\frac{2k}{k+m_1-1}}\\ 
&\leq \const{gn23} \|\nabla(u+1)^{\frac{k+m_1-1}{2}}\|_{L^2(\Omega)}^{\sigma_1\theta_1} \|(u+1)^{\frac{k+m_1-1}{2}}\|_{L^\frac{2}{k+m_1-1}(\Omega)}^{\sigma_1(1-\theta_1)} + \const{gn23} \|(u+1)^{\frac{k+m_1-1}{2}}\|_{L^\frac{2}{k+m_1-1}(\Omega)}^{\sigma_1}\\
&\leq \const{gn24} \left(\int_\Omega |\nabla (u+1)^{\frac{k+m_1-1}{2}}|^2 \right)^\frac{k\theta_1}{k+m_1-1}+\const{gn25}
\leq \const{gn26} \left(\int_\Omega |\nabla (u+1)^{\frac{k+m_1-1}{2}}|^2 +1\right)^{\frac{k\theta_1}{k+m_1-1}},
\end{split}
\end{equation*}
from which we deduce the claim after basic manipulations.

Now, we focus on the first bound of \eqref{GN3} and \eqref{GN4}, respectively. Due to $\delta > \frac{n \alpha}{n+1}$, by exploiting \eqref{theta2}, \eqref{sigmatheta2} and the boundedness of the mass 
in \eqref{MassBounded_grad}, a further application of the Gagliardo--Nirenberg and Young's inequalities leads to
\begin{equation}\label{SupportoPerCasoLimite}
\begin{split}
\hat{c}\int_\Omega (u+1)^{k+\alpha-1}&=\hat{c}\|(u+1)^{\frac{k-1+\delta}{\delta}}\|_{L^{\frac{\delta(k+\alpha-1)}{k-1+\delta}}(\Omega)}^{\frac{\delta(k+\alpha-1)}{k-1+\delta}}\\ 
&\leq \const{gn32} \|\nabla(u+1)^{\frac{k-1+\delta}{\delta}}\|_{L^\delta(\Omega)}^{\sigma_2\theta_2} \|(u+1)^{\frac{k-1+\delta}{\delta}}\|_{L^\frac{\delta}{k-1+\delta}(\Omega)}^{\sigma_2(1-\theta_2)}
+ \const{gn32} \|(u+1)^{\frac{k-1+\delta}{\delta}}\|_{L^\frac{\delta}{k-1+\delta}(\Omega)}^{\sigma_2}\\
&\leq \const{gn33} \left(\int_\Omega |\nabla (u+1)^{\frac{k-1+\delta}{\delta}}|^\delta \right)^\frac{\sigma_2\theta_2}{\delta}+\const{gn33}\\
&\leq \frac{kc}{2} \left(\frac{k-1+\delta}{\delta}\right)^{-\delta} \int_\Omega |\nabla (u+1)^{\frac{k-1+\delta}{\delta}}|^\delta+\const{gn3}
 \quad \text{on }(0,\TM).
\end{split}
\end{equation}
In a similarly way as before, for $\delta >\frac{n(m_2+l)}{n+1}$, by applying \eqref{theta3}, \eqref{sigmatheta3} and the boundedness of the mass in \eqref{MassBounded_grad}, we can estimate the integral term $\hat{c} \int_\Omega (u+1)^{k+m_2+l-1}$ on $(0,\TM)$ as
\begin{equation*}
\begin{split}
\hat{c} \int_\Omega (u+1)^{k+m_2+l-1}&=\hat{c} \|(u+1)^{\frac{k-1+\delta}{\delta}}\|_{L^{\frac{\delta(k+m_2+l-1)}{k-1+\delta}}(\Omega)}^{\frac{\delta(k+m_2+l-1)}{k-1+\delta}}\\ 
&\leq \const{gn42} \|\nabla(u+1)^{\frac{k-1+\delta}{\delta}}\|_{L^\delta(\Omega)}^{\sigma_3\theta_3} \|(u+1)^{\frac{k-1+\delta}{\delta}}\|_{L^\frac{\delta}{k-1+\delta}(\Omega)}^{\sigma_3(1-\theta_3)}
+ \const{gn42} \|(u+1)^{\frac{k-1+\delta}{\delta}}\|_{L^\frac{\delta}{k-1+\delta}(\Omega)}^{\sigma_3}\\
&\leq \const{gn43} \left(\int_\Omega |\nabla (u+1)^{\frac{k-1+\delta}{\delta}}|^\delta \right)^\frac{\sigma_3 \theta_3}{\delta}+\const{gn44}
\leq  \frac{kc}{2} \left(\frac{k-1+\delta}{\delta}\right)^{-\delta} \int_\Omega |\nabla (u+1)^{\frac{k-1+\delta}{\delta}}|^\delta+\const{gn4}.
 \end{split}
\end{equation*}
As to the limit cases of relations \eqref{GN3} and \eqref{GN4}, we focus only to the first, being similar the other. With the inequality \eqref{SupportoPerCasoLimite} in our hands, we essentially observe that for $\delta=\frac{n\alpha}{n+1}$ the exponent $\frac{\sigma_2\theta_2}{\delta}$ equals $1$, and henceforth the Young inequality in the last step is no more applicable; as a consequence we have
\begin{equation*}
\begin{split}
\hat{c}\int_\Omega (u+1)^{k+\alpha-1}&=\hat{c}\|(u+1)^{\frac{k-1+\delta}{\delta}}\|_{L^{\frac{\delta(k+\alpha-1)}{k-1+\delta}}(\Omega)}^{\frac{\delta(k+\alpha-1)}{k-1+\delta}}\\ 
%\leq& \const{gn32} \|\nabla(u+1)^{\frac{k-1+\delta}{\delta}}\|_{L^\delta(\Omega)}^{\sigma_2\theta_2} %\|(u+1)^{\frac{k-1+\delta}{\delta}}\|_{L^\frac{\delta}{k-1+\delta}(\Omega)}^{\sigma_2(1-\theta_2)}
%+ \hat{c} \|(u+1)^{\frac{k-1+\delta}{\delta}}\|_{L^\frac{\delta}{k-1+\delta}(\Omega)}^{\sigma_2}%\\
&\leq \hat{c}(2 C_{GN})^{\sigma_2}\left(\int_\Omega |\nabla (u+1)^{\frac{k-1+\delta}{\delta}}|^\delta \right)^\frac{\sigma_2\theta_2}{\delta}\left(\int_\Omega (u+1)\right)^{(k+\alpha-1)(1-\theta_2)}+\const{gn33}\\
&\leq  \hat{c}\, (2 C_{GN})^{\sigma_2} (M_1+|\Omega|)^{\frac{\alpha}{n+1}} \int_\Omega |\nabla (u+1)^{\frac{k-1+\delta}{\delta}}|^\delta+\const{gn03}
 \quad \text{on }(0,\TM).
\end{split}
\end{equation*}
\end{proof}
\end{lemma}
In the next result we will, indeed, rely on Lemma \ref{GagliardoIneqLemma}.
\begin{lemma}\label{GNandInterpBisLemma}
Let the assumptions of Lemma \ref{ThetaSigma} be valid and let $k_*>1$ be the value described in Convention \ref{Convenzione}. Then for every $\hat{c}>0$, $k>k_*$ and all $t\in (0,\TM)$ we have for $\epsilon>0$
\begin{equation}\label{Interpolation1}
\begin{split}
\hat{c} \int_\Omega (u+1)^{k+m_2+l-1} &\le 
\begin{cases}
\displaystyle \frac{k(k-1)}{(k+m_1-1)^2}
\displaystyle \int_\Omega |\nabla (u+1)^{\frac{k+m_1-1}{2}}|^2 +\epsilon
\int_\Omega (u+1)^{\delta} \int_\Omega (u+1)^{k-1}+\const{i1} \\
\quad \quad \quad \quad \quad \quad\quad \quad \quad \quad \quad \quad \textrm{ if } 
\displaystyle \delta\geq m_1 \textrm{ and } \displaystyle  \delta > \frac{n(m_2-m_1+l)+2(m_2+l)}{2},\\ \\
\displaystyle \frac{k(k-1)}{(k+m_1-1)^2}
\int_\Omega |\nabla (u+1)^{\frac{k+m_1-1}{2}}|^2 +\hat{c}C_0
\int_\Omega (u+1)^{\delta} \int_\Omega (u+1)^{k-1}+\const{0i1} \\
\quad \quad \quad \quad \quad \quad \quad \quad \quad\quad \quad \quad   \textrm{ if } 
\displaystyle\delta\geq m_1 \textrm{ and } \displaystyle \delta = \frac{n(m_2-m_1+l)+2(m_2+l)}{2},
\end{cases}
\end{split}
\end{equation}
\begin{equation}\label{InterpolationP1}
\begin{split}
\hat{c} \int_\Omega (u+1)^{k+m_2+l-1} &\le 
\begin{cases}
\displaystyle \frac{k(k-1)}{(k+m_1-1)^2}
\displaystyle \int_\Omega |\nabla (u+1)^{\frac{k+m_1-1}{2}}|^2 +\epsilon
\int_\Omega (u+1)^{\beta} \int_\Omega (u+1)^{k+\alpha-1}+\const{ip1} \\
\quad \quad \quad \quad \quad \quad\quad \quad \quad \quad  \textrm{ if } 
\displaystyle \beta\geq m_1-\alpha \textrm{ and } \beta > \frac{n(m_2-m_1+l)+2(m_2+l-\alpha)}{2},\\ \\
\displaystyle \frac{k(k-1)}{(k+m_1-1)^2}\int_\Omega |\nabla (u+1)^{\frac{k+m_1-1}{2}}|^2 +\hat{c} C_0
\int_\Omega (u+1)^{\beta} \int_\Omega (u+1)^{k+\alpha-1}+\const{0ip1} \\
\quad \quad \quad \quad \quad \quad\quad \quad \quad \quad  \textrm{ if } 
\displaystyle \beta\geq m_1-\alpha \textrm{ and } \beta = \frac{n(m_2-m_1+l)+2(m_2+l-\alpha)}{2},
\end{cases}
\end{split}
\end{equation}
\begin{equation}\label{InterpolationP2}
\begin{split}
\hat{c} \int_\Omega (u+1)^{k+\alpha-1} &\le 
\begin{cases}
\displaystyle \frac{k(k-1)}{(k+m_1-1)^2}
\displaystyle \int_\Omega |\nabla (u+1)^{\frac{k+m_1-1}{2}}|^2 +\epsilon
\int_\Omega (u+1)^{\beta} \int_\Omega (u+1)^{k+\alpha-1}+\const{ip2} \\
\quad \quad \quad \quad \quad \quad\quad \quad \quad \quad \quad \quad \textrm{ if } 
\displaystyle \beta\geq m_1-\alpha \textrm{ and } \beta > \frac{n(\alpha-m_1)}{2},\\ \\
\displaystyle \frac{k(k-1)}{(k+m_1-1)^2}
\int_\Omega |\nabla (u+1)^{\frac{k+m_1-1}{2}}|^2 +\hat{c} C_0
\int_\Omega (u+1)^{\beta} \int_\Omega (u+1)^{k+\alpha-1}+\const{0ip2} \\
\quad \quad \quad \quad \quad \quad\quad \quad \quad \quad \quad \quad \textrm{ if } 
\displaystyle \beta\geq m_1-\alpha \textrm{ and } \beta = \frac{n(\alpha-m_1)}{2}.
\end{cases}
\end{split}
\end{equation}
\begin{proof}
In order to show the non-limit case of relations \eqref{Interpolation1},  
\eqref{InterpolationP1} and \eqref{InterpolationP2}, we will  
exploit Lemma \ref{GagliardoIneqLemma} with $\psi:=(u+1)^{\frac{k+m_1-1}{2}}$ and proper $q$ and $r$, exactly fulfilling the restrictions in the lemma itself. 

First we focus our attention on the proof of \eqref{Interpolation1}. By virtue of the assumption $\delta \geq m_1$, for any $k>k_*$ it is possible to set 
\begin{equation}\label{kPrimoInterp1}
k':=\frac{k+\delta-1}{2},
\end{equation}
which satisfies 
\begin{equation*}%\label{k'condiInterp1}
\max\left\{\delta,\ \frac{k+m_1-1}{2},\ \frac{n(m_2+l-m_1)}{2}\right\}<k'<\min\{k-1,\ k+m_2+l-1\}.
\end{equation*}
In this way, 
\begin{equation}\label{qr1}
q:=\frac{2(k+m_2+l-1)}{k+m_1-1},\ r:=\frac{2k'}{k+m_1-1}
\end{equation}
are consistent choices.
%a number of calculations yields 
%\begin{equation}\label{Cond_pqr}
%\begin{split}
%&1\le r<q<p \quad \text{and} \quad \frac{q}{r}<\frac{2}{r}+1-\frac{2}{p} \quad \text{if} \quad n \geq %3,\\
%&1\le r<q<p \quad \text{and} \quad \frac{q}{r}<\frac{2}{r}+2 \quad \text{if} \quad n=1,\\
%&1\le r<q<p \quad \text{and} \quad \frac{q}{r}<\frac{2}{r}+1 \quad \text{if} \quad n=2.
%\end{split}
%\end{equation}
Therefore we infer from \eqref{Bian} that 
\begin{equation}\label{u^k+m_2+l-1Interp1}
\begin{split}
\hat{c} \int_\Omega (u+1)^{k+m_2+l-1}&= \hat{c}\|(u+1)^{\frac{k+m_1-1}{2}}\|_{L^{\frac{2(k+m_2+l-1)}{k+m_1-1}}(\Omega)}^{\frac{2(k+m_2+l-1)}{k+m_1-1}}
\leq \frac{k(k-1)}{2(k+m_1-1)^2}  \int_\Omega |\nabla (u+1)^{\frac{k+m_1-1}{2}}|^2\\
& \quad+ \hat{c} \int_\Omega (u+1)^{k+m_1-1}+ \hat{c} C_0 \left(\int_\Omega (u+1)^{k'}\right)^{\frac{\gamma}{r}}\quad  \mbox{for\ all}\ t \in (0, \TM). 
\end{split}
\end{equation}
Here, the interpolation inequality (see \cite[page 93]{BrezisBook}) leads to 
\begin{equation}\label{u^kpInterp1}
\begin{split}
&\hat{c} C_0 \left(\int_\Omega (u+1)^{k'}\right)^{\frac{\gamma}{r}}=\hat{c} C_0\|(u+1)\|_{L^{k'}(\Omega)}^{b_1}
\le \hat{c} C_0 \|(u+1)\|_{L^\delta(\Omega)}^{a_1 b_1} \|(u+1)\|_{L^{k-1}(\Omega)}^{(1-a_1)b_1}\\
&=\hat{c} C_0 \left(\|(u+1)\|_{L^\delta(\Omega)}^\delta \|(u+1)\|_{L^{k-1}(\Omega)}^{k-1}\right)^{\frac{a_1 b_1}{\delta}}\|(u+1)\|_{L^{k-1}(\Omega)}
^{\left[1-a_1-\frac{a_1(k-1)}{\delta}\right]b_1} \quad \textrm{on } (0,\TM),
\end{split}
\end{equation}
where 
\begin{equation*}%\label{DefinitionB1A1}
b_1=b_1(q):=\frac{k'\gamma(q)}{r}=\frac{k'\gamma}{r}, \quad a_1:=\frac{\frac{1}{k'}-\frac{1}{k-1}}{\frac{1}{\delta}-\frac{1}{k-1}}\in (0,1).
\end{equation*}
We note that recalling the expression of $k'$ in \eqref{kPrimoInterp1} and the constrain $\delta >\frac{n(m_2-m_1+l)+2(m_2+l)}{2}$, some computations provide 
\begin{equation*}%\label{pow_zero}
\left[1-a_1-\frac{a_1(k-1)}{\delta}\right]b_1=0 \quad \textrm{and} \quad 
\frac{a_1 b_1}{\delta}<1. 
\end{equation*}
As a consequence, we can invoke Young's inequality so that  
relation \eqref{u^kpInterp1} reads 
\begin{equation}\label{SupportoDue}
\begin{split}
\hat{c} C_0\left(\int_\Omega (u+1)^{k'}\right)^{\frac{\gamma}{r}} &\le \hat{c} C_0\left(\|(u+1)\|_{L^\delta(\Omega)}^\delta \|(u+1)\|_{L^{k-1}(\Omega)}^{k-1}\right)
^{\frac{a_1 b_1}{\delta}}\\
&\le \epsilon \int_\Omega (u+1)^{\delta} \int_\Omega (u+1)^{k-1} +\const{j} \quad  \mbox{for\ all}\ t \in (0, \TM),
\end{split}
\end{equation}
which in conjunction with \eqref{u^k+m_2+l-1Interp1} implies for all $t \in (0, \TM)$, 
\begin{equation}\label{u^k+m_2+l-1_2Interp1}
\begin{split}
\hat{c} \int_\Omega (u+1)^{k+m_2+l-1} &\le \frac{k(k-1)}{2(k+m_1-1)^2}  \int_\Omega |\nabla (u+1)^{\frac{k+m_1-1}{2}}|^2+ \hat{c} \int_\Omega (u+1)^{k+m_1-1}\\
& \quad + \epsilon \int_\Omega (u+1)^{\delta} \int_\Omega (u+1)^{k-1} +\const{j}.
\end{split}
\end{equation}
Now, we can treat the second integral term of the right-hand side of \eqref{u^k+m_2+l-1_2Interp1} by employing the Gagliardo--Nirenberg and Young's inequalities (recalling \eqref{theta4} and the boundedness of the mass in Lemma \ref{localSol} for problems \eqref{problem_int} and \eqref{problem_grad}), entailing
\begin{equation*}%\label{TermInterp1}
\begin{split}
&\hat{c} \int_\Omega (u+1)^{k+m_1-1}= \hat{c} \|(u+1)^{\frac{k+m_1-1}{2}}\|_{L^2(\Omega)}^2\\ 
\leq& \const{l} \|\nabla(u+1)^{\frac{k+m_1-1}{2}}\|_{L^2(\Omega)}^{2\theta_4}
\|(u+1)^{\frac{k+m_1-1}{2}}\|_{L^\frac{2}{k+m_1-1}(\Omega)}^{2(1-\theta_4)}
 + \const{l} \|(u+1)^{\frac{k+m_1-1}{2}}\|_{L^\frac{2}{k+m_1-1}(\Omega)}^2\\
 \leq& \const{m} \left(\int_\Omega |\nabla (u+1)^{\frac{k+m_1-1}{2}}|^2 \right)^{\theta_4}+\const{n} \leq \frac{k(k-1)}{2(k+m_1-1)^2} \int_\Omega |\nabla (u+1)^{\frac{k+m_1-1}{2}}|^2 +\const{n}
 \quad \text{on }(0,\TM),
\end{split}
\end{equation*}
which together with \eqref{u^k+m_2+l-1_2Interp1} gives the claim.

In turn, let us show estimate \eqref{InterpolationP1}. From the condition $\beta \geq m_1-\alpha$, we have that for any $k>k_*$ 
\begin{equation}\label{kPrimoInterpP1}
k':=\frac{k+\alpha+\beta-1}{2}
\end{equation}
satisfies 
\begin{equation*}%\label{k'condiInterpP1}
\max\left\{\beta,\ \frac{k+m_1-1}{2},\ \frac{n(m_2+l-m_1)}{2}\right\}<k'<
\min\{k+\alpha-1,\ k+m_2+l-1\}.
\end{equation*}
In this way, for $q$ and $r$ taken accordingly to \eqref{qr1}, by reasoning as in the proof of the first inequality in \eqref{Interpolation1}, we derive for all $t \in (0, \TM)$ 
\begin{equation}\label{u^k+m_2+l-1Interp1P}
\begin{split}
\hat{c} \int_\Omega (u+1)^{k+m_2+l-1}&= \hat{c}\|(u+1)^{\frac{k+m_1-1}{2}}\|_{L^{\frac{2(k+m_2+l-1)}{k+m_1-1}}(\Omega)}^{\frac{2(k+m_2+l-1)}{k+m_1-1}}\\
&\le  \frac{k(k-1)}{(k+m_1-1)^2} \int_\Omega |\nabla (u+1)^{\frac{k+m_1-1}{2}}|^2
+ \hat{c} C_0 \left(\int_\Omega (u+1)^{k'}\right)^{\frac{\gamma}{r}}\\
& \leq  \frac{k(k-1)}{(k+m_1-1)^2} \int_\Omega |\nabla (u+1)^{\frac{k+m_1-1}{2}}|^2\\
& \quad + \hat{c} C_0 \left(\|(u+1)\|_{L^\beta(\Omega)}^\beta \|(u+1)\|_{L^{k+\alpha-1}(\Omega)}^{k+\alpha-1}\right)^{\frac{a_2 b_2}{\beta}}\|(u+1)\|_{L^{k+\alpha-1}(\Omega)}^{\left[1-a_2-\frac{a_2(k+\alpha-1)}{\beta}\right]b_2},
 \end{split}
\end{equation}
where 
\begin{equation*}%\label{DefinitionB3A3}
b_2=b_2(q):=\frac{k'\gamma(q)}{r}=\frac{k'\gamma}{r}, \quad a_2:=\frac{\frac{1}{k'}-\frac{1}{k+\alpha-1}}{\frac{1}{\beta}-\frac{1}{k+\alpha-1}}\in (0,1).
\end{equation*}
Additionally, taking in mind the expression of $k'$ given in \eqref{kPrimoInterpP1} and the restriction $\beta> \frac{n(m_2-m_1+l)+2(m_2+l-\alpha)}{2}$, we can establish that
\begin{equation*}%\label{pow_zero}
\left[1-a_2-\frac{a_2(k+\alpha-1)}{\beta}\right]b_2=0 \quad \textrm{and} \quad 
\frac{a_2 b_2}{\beta}<1. 
\end{equation*}
Again an application of Young's inequality allows to reduce relation \eqref{u^k+m_2+l-1Interp1P} for all $t \in (0, \TM)$ into
 \begin{equation*}%\label{u^k+m_2+l-1_2Interp1P}
 \hat{c} \int_\Omega (u+1)^{k+m_2+l-1} \le \frac{k(k-1)}{(k+m_1-1)^2} \int_\Omega |\nabla (u+1)^{\frac{k+m_1-1}{2}}|^2+\epsilon \int_\Omega (u+1)^{\beta} \int_\Omega (u+1)^{k+\alpha-1}+\const{ip1}.
\end{equation*}
Finally, let us deduce estimate \eqref{InterpolationP2}. From the assumption on 
$\beta$, for any $k>k_*$ the value $k'$ in \eqref{kPrimoInterpP1}
complies with
\begin{equation*}%\label{k'condiInterp2P}
\max\left\{\beta,\ \frac{k+m_1-1}{2}\right\}<k'<k+\alpha-1.
\end{equation*}
In this way, for  
\begin{equation*}%\label{qr2}
q:=\frac{2(k+\alpha-1)}{k+m_1-1},\ r:=\frac{2k'}{k+m_1-1},
\end{equation*}
we derive for all $ t \in (0, \TM)$
\begin{equation}\label{u^k+alpha-1Interp2P}
\begin{split}
\hat{c} \int_\Omega (u+1)^{k+\alpha-1}&= \hat{c} \|(u+1)^{\frac{k+m_1-1}{2}}\|_{L^{\frac{2(k+\alpha-1)}{k+m_1-1}}(\Omega)}^{\frac{2(k+\alpha-1)}{k+m_1-1}}\\
&\le \frac{k(k-1)}{(k+m_1-1)^2} \int_\Omega |\nabla (u+1)^{\frac{k+m_1-1}{2}}|^2
+\hat{c} C_0 \left(\int_\Omega (u+1)^{k'}\right)^{\frac{\gamma}{r}}\\
& \le  \frac{k(k-1)}{(k+m_1-1)^2} \int_\Omega |\nabla (u+1)^{\frac{k+m_1-1}{2}}|^2\\
&\quad  + \hat{c} C_0 \left(\|(u+1)\|_{L^\beta(\Omega)}^\beta \|(u+1)\|_{L^{k+\alpha-1}(\Omega)}^{k+\alpha-1}\right)^{\frac{a_3 b_3}{\beta}}\|(u+1)\|_{L^{k+\alpha-1}(\Omega)}
^{\left[1-a_3-\frac{a_3(k+\alpha-1)}{\beta}\right]b_3},
\end{split}
\end{equation}
where 
\begin{equation*}%\label{DefinitionB4A4}
b_3=b_3(q):=\frac{k'\gamma(q)}{r}=\frac{k'\gamma}{r}, \quad a_3:=a_2\in (0,1) \quad (a_2  \textrm{ has been defined above}).
\end{equation*}
In a very similar way, taking into account $k'$ in \eqref{kPrimoInterpP1} and the assumption on $\beta$, we have 
\begin{equation*}%\label{pow_zero}
\left[1-a_3-\frac{a_3(k+\alpha-1)}{\beta}\right]b_3=0 \quad \textrm{and} \quad 
\frac{a_3 b_3}{\beta}<1,
\end{equation*}
and Young's inequality makes that  relation \eqref{u^k+alpha-1Interp2P} is turned for all $t \in (0, \TM)$ into
\begin{equation*}%\label{u^k+m_2+l-1_2Interp1P}
\hat{c} \int_\Omega (u+1)^{k+\alpha-1} \le  \frac{k(k-1)}{(k+m_1-1)^2} \int_\Omega |\nabla (u+1)^{\frac{k+m_1-1}{2}}|^2+\epsilon \int_\Omega (u+1)^{\beta} \int_\Omega (u+1)^{k+\alpha-1}+\const{ip2}.
\end{equation*}
Let us show only the limit case for \eqref{Interpolation1}. The sole modification with respect the non-limit situation regards the analysis of the term \eqref{SupportoDue}. We have that $\frac{a_1 b_1}{\delta}=1$, in view of $\delta =\frac{n(m_2-m_1+l)+2(m_2+l)}{2}$; henceforth \eqref{SupportoDue} is reorganized as
\begin{equation*}
\begin{split}
\hat{c} C_0\left(\int_\Omega (u+1)^{k'}\right)^{\frac{\gamma}{r}} 
%\hat{c} C_0\left(\|(u+1)\|_{L^\delta(\Omega)}^\delta \|(u+1)\|_{L^{k-1}(\Omega)}^{k-1}\right)
%^{\frac{a_1 b_1}{\delta}}\\
\le \hat{c} C_0\int_\Omega (u+1)^{\delta} \int_\Omega (u+1)^{k-1} \quad  \mbox{for\ all}\ t \in (0, \TM).
\end{split}
\end{equation*}
\end{proof}
\end{lemma}
\subsection{A priori boundedness of $u$ for model \eqref{problem_int}}
In order to investigate further properties of the solution $(u,v)$ to the models of our interest, let us use that $u$ and $v$ actually solve the related equations.
\begin{lemma}\label{StimaGenerale}
Let the assumptions of Lemma \ref{ThetaSigma} be valid and let $k_*>1$ be  the value described in Convention \ref{Convenzione}. Then, for $k>k_*$, the functional $y(t)$ defined in \eqref{Func} fulfills  for all $\epsilon>0$ and for $F_{m_2}(u)$ as in \eqref{Def_Fm2}
\begin{equation}\label{Stima1}
\begin{split}
y'(t)&\leq -\frac{4k(k-1)}{(k+m_1-1)^2} \int_\Omega |\nabla (u+1)^{\frac{k+m_1-1}{2}}|^2 - k(k-1) \chi \int_\Omega F_{m_2}(u)\Delta v + ka \int_\Omega (u+1)^{k+\alpha-1}\\
& \quad +\left(\epsilon-\frac{kb}{2^{\beta-1}}\right) \int_\Omega (u+1)^{k+\beta-1}
- \frac{kc}{2^{\delta-1}} \int_\Omega (u+1)^{k-1} \int_\Omega (u+1)^{\delta} + \const{s1} \quad \mbox{on}\ (0, \TM).
\end{split}
\end{equation}
\begin{proof}
 For all $k>k_*$, we have from the first equation in \eqref{problem_int} and integration by parts that 
\begin{equation*}
\begin{split}
y'(t)=\frac{d}{dt}\int_\Omega (u+1)^k &=k\int_\Omega (u+1)^{k-1} u_t = -k(k-1)\int_\Omega (u+1)^{k+m_1-3} |\nabla u|^2\\ 
&\quad + k(k-1) \chi \int_\Omega u(u+1)^{k+m_2-3} \nabla u \cdot \nabla v +ka \int_\Omega u^\alpha (u+1)^{k-1} - kb  \int_\Omega u^\beta (u+1)^{k-1}\\ 
&\quad - kc \int_\Omega (u+1)^{k-1} \int_\Omega u^{\delta}
\quad  \mbox{for all }\ t \in (0, \TM). 
\end{split}
\end{equation*}
Taking into account the definition and the bound of $F_{m_2}(u)$ in \eqref{Def_Fm2} and \eqref{Fm2} and applying the divergence theorem, the previous inequality can be rewritten as 
\begin{equation}\label{Stima2}
\begin{split}
y'(t) &= -\frac{4k(k-1)}{(k+m_1-1)^2} \int_\Omega |\nabla (u+1)^{\frac{k+m_1-1}{2}}|^2 
- k(k-1) \chi \int_\Omega F_{m_2}(u) \Delta v \\
&\quad +ka \int_\Omega u^\alpha (u+1)^{k-1}- kb  \int_\Omega u^\beta (u+1)^{k-1} - kc \int_\Omega (u+1)^{k-1} \int_\Omega u^{\delta} \quad \textrm{on } (0, \TM). 
\end{split}
\end{equation}
By exploiting the trivial relation $u \leq u+1$ and inequality \eqref{ABIn}, estimate \eqref{Stima2} becomes 
\begin{equation}\label{Stima3}
\begin{split}
y'(t) &\leq -\frac{4k(k-1)}{(k+m_1-1)^2} \int_\Omega |\nabla (u+1)^{\frac{k+m_1-1}{2}}|^2 - k(k-1) \chi \int_\Omega F_{m_2}(u) \Delta v
\\
& \quad + ka \int_\Omega (u+1)^{k+\alpha-1} + (kb+kc|\Omega|) \int_\Omega (u+1)^{k-1} 
-\frac{kb}{2^{\beta-1}}  \int_\Omega (u+1)^{k+\beta-1}\\ 
& \quad - \frac{kc}{2^{\delta-1}} \int_\Omega (u+1)^{k-1} \int_\Omega (u+1)^{\delta} \quad \mbox{for all } t\in (0, \TM).
\end{split}
\end{equation}
Due to $\beta \geq 1$, the Young inequality gives 
\begin{equation}\label{Y1}
(kb+kc|\Omega|) \int_\Omega (u+1)^{k-1} \leq \epsilon \int_\Omega (u+1)^{k+\beta-1} + \const{s1} \quad \textrm{on } (0,\TM).
\end{equation}
By inserting relation \eqref{Y1} into estimate \eqref{Stima3}, we conclude the proof. 
\end{proof}
\end{lemma}

\begin{lemma}\label{LemmaProblemIntEll1} 
Let the assumptions of Lemma \ref{ThetaSigma} be valid, let $k_*>1$ be the value described in Convention \ref{Convenzione} and let $k>k_*$. Additionally, let $f$ comply with 
\eqref{Def_f} and let $\tau \in \{0,1\}$, $c>0$, $\beta>\alpha\geq 1$ and $\delta\geq \max\{1,m_1\}$ satisfy either
%\begin{itemize}
%\item if $c=0$ the conditions 
%\begin{equation}\label{czero}
%\beta > m_2+l  \quad \textrm{or} \quad m_2+l<m_1+ \frac{2}{n};
%\end{equation}
%\item if $c>0$
 % \item 
 \begin{equation}\label{cdivzero1}
 \delta >\frac{n(m_2-m_1+l)+2(m_2+l)}{2}, 
\end{equation}
or 
%\item 
\begin{equation}\label{LimitInt1}
 \delta =\frac{n(m_2-m_1+l)+2(m_2+l)}{2} \quad \textrm{and} \quad c \geq \frac{2^{\delta-1}(k-1)\chi K C_0}{k+m_2-1}\left(\tau C_P^{\frac{l}{k+m_2+l-1}}+(1-\tau)\right).
 \end{equation}
%\end{itemize}
Then $u\in L^\infty((0,\TM);L^k(\Omega))$ for all $k>k_*$. 
%\begin{equation*}
%\int_\Omega u^k \le L_0 \quad  \mbox{for\ all}\ t \in (0, \TM). %
%\end{equation*}
\begin{proof}
We distinguish the cases $\tau=0$ and $\tau=1$.
\subsubsection*{$\triangleright$ $\tau=0$}
Taking into account the definition of $f$ in \eqref{Def_f}, the bound of $F_{m_2}$ in \eqref{Fm2} and exploiting the second equation in \eqref{problem_int} , we have on $(0,\TM)$ that for all $k>k_*$
\begin{equation}\label{StimaFm2}
- k(k-1)\chi \int _\Omega F_{m_2}(u)\Delta v=- k(k-1)\chi \int _\Omega F_{m_2}(u)(v-f(u))\leq \frac{k(k-1)\chi K}{k+m_2-1} \int _\Omega (u+1)^{k+m_2+l-1},
\end{equation}
where we have also neglected nonpositive terms. 
%\begin{equation}\label{S1}
%\begin{split}
%y'(t) &\leq -\frac{4k(k-1)}{(k+m_1-1)^2} \int_\Omega |\nabla (u+1)^{\frac{k+m_1-1}{2}}|^2 - k(k-1)\chi \int _\Omega F_{m_2}(u)(v-f(u))\\
%& \; + ka \int_\Omega (u+1)^{k+\alpha-1} -\frac{kb}{2^{\beta}}  \int_\Omega (u+1)^{k+\beta-1}- \frac{kc}{2^{\delta-1}} \int_\Omega (u+1)^{k-1} \int_\Omega (u+1)^{\delta}\\
%&\leq -\frac{4k(k-1)}{(k+m_1-1)^2} \int_\Omega |\nabla (u+1)^{\frac{k+m_1-1}{2}}|^2 
%+ \frac{k(k-1)\chi K}{k+m_2-1} \int _\Omega (u+1)^{k+m_2+l-1}\\
%& \; + ka \int_\Omega (u+1)^{k+\alpha-1} -\frac{kb}{2^{\beta}} \int_\Omega (u+1)^{k+\beta-1}- %\frac{kc}{2^{\delta-1}} \int_\Omega (u+1)^{k-1} \int_\Omega (u+1)^{\delta}. 
% \end{split}
%\end{equation}
On the other hand, from $\beta >\alpha$ one can apply the Young inequality so entailing
\begin{equation}\label{Y2}
ka \int_\Omega (u+1)^{k+\alpha-1} \leq \epsilon \int_\Omega (u+1)^{k+\beta-1} + \const{b} \quad \textrm{on } (0,\TM).
\end{equation}
These two gained relations, plugged into \eqref{Stima1}, imply 
\begin{equation}\label{Sf1}
\begin{split}
y'(t) &\leq -\frac{4k(k-1)}{(k+m_1-1)^2} \int_\Omega |\nabla (u+1)^{\frac{k+m_1-1}{2}}|^2 
+ \frac{k(k-1)\chi K}{k+m_2-1} \int _\Omega (u+1)^{k+m_2+l-1}\\
& \quad +\left(\epsilon-\frac{kb}{2^{\beta-1}}\right) \int_\Omega (u+1)^{k+\beta-1}- \frac{kc}{2^{\delta-1}} \int_\Omega (u+1)^{k-1} \int_\Omega (u+1)^{\delta} + \const{Bb} \quad \textrm{for all } t \in (0,\TM). 
 \end{split}
\end{equation}
 Assumption \eqref{cdivzero1} allows us to apply the first inequality in \eqref{Interpolation1}, which in conjunction with \eqref{Sf1} and for suitable $\epsilon$ gives the inequality 
 \begin{equation}\label{uk_estimate4}
y'(t) \leq - \const{h} \int_\Omega |\nabla (u+1)^{\frac{k+m_1-1}{2}}|^2 + \const{si}  \quad \mbox{on}\ (0, \TM).
\end{equation}
Therefore, by virtue of relation \eqref{GN2}, the above inequality leads to this initial problem 
\begin{equation*}%\label{eq:sistemafinale}
%\begin{dcases}
y'(t) \leq \const{ac1} - \const{ac2} y(t)^{\frac{k+m_1-1}{k\theta_1}} \quad \text{on }(0,\TM),\quad
y(0)=\int_\Omega (u_0+1)^p,
%\end{dcases}
\end{equation*}
so that an application of Lemma \ref{LemmaODI-Comparison} implies $\sup_{t\in (0,\TM)}\int_\Omega u^k \leq y(t)<\infty$.

To deal with the limit situation, by considering inequality \eqref{Sf1} and, again, by neglecting the nonpositive term associated to the logistic, we obtain for all $t \in (0,\TM)$
\begin{equation*}\label{Ctau0}
\begin{split}
y'(t) + \frac{kc}{2^{\delta-1}} \int_\Omega (u+1)^{k-1} \int_\Omega (u+1)^{\delta} &\leq -\frac{4k(k-1)}{(k+m_1-1)^2} \int_\Omega |\nabla (u+1)^{\frac{k+m_1-1}{2}}|^2 \\
& \quad 
+ \frac{k(k-1)\chi K}{k+m_2-1} \int _\Omega (u+1)^{k+m_2+l-1}
+ \const{Bb}. 
 \end{split}
\end{equation*}
Since $\delta = \frac{n(m_2-m_1+l)+2(m_2+l)}{2}$, we plug in the above the second estimate of \eqref{Interpolation1} with the particular choice $\hat{c}=\frac{k(k-1)\chi K}{k+m_2-1}$, so concluding by virtue of the largeness assumption $c \geq \frac{2^{\delta-1}(k-1)\chi K C_0}{k+m_2-1}$.
\subsubsection*{$\triangleright$ $\tau=1$}
First of all, we focus on the integral term $\int_\Omega F_{m_2}(u) \Delta v$, which can be estimated by exploiting Young's inequality due to $l>0$ and bound on $F_{m_2}$ in \eqref{Fm2}, entailing 
\begin{equation}\label{ParabSt1}
\begin{split}
&-k(k-1)\chi \int_\Omega F_{m_2}(u) \Delta v \leq k(k-1)\chi \int_\Omega F_{m_2}(u) 
|\Delta v| \leq \frac{k(k-1)\chi}{k+m_2-1} \int_\Omega (u+1)^{k+m_2-1} |\Delta v|\\
&\leq \frac{k(k-1)\chi K C_P^{\frac{l}{k+m_2+l-1}}}{k+m_2+l-1} 
\int_\Omega (u+1)^{k+m_2+l-1}\\ 
& \qquad 
+ \frac{k(k-1)\chi l}{(k+m_2-1)(k+m_2+l-1)}  \left(K C_P^{\frac{l}{k+m_2+l-1}}\right)^{-\frac{k+m_2-1}{l}} \int_\Omega |\Delta v|^{\frac{k+m_2+l-1}{l}} \quad \textrm{for all } t \in (0,\TM). 
\end{split}
\end{equation}
At this point, by applying Lemma \ref{lem:MaxReg} to the second equation of \eqref{problem_int} with $\psi=v$, $q=\frac{k+m_2+l-1}{l}$ and $h=f$ given in \eqref{Def_f}, we deduce that 
\begin{equation}\label{ParabSt2} 
\begin{split}
&\frac{k(k-1)\chi l}{(k+m_2-1)(k+m_2+l-1)}  \left(K C_P^{\frac{l}{k+m_2+l-1}}\right)^{-\frac{k+m_2-1}{l}} \int_0^t e^s \left(\int_\Omega |\Delta v(\cdot,s)|^{\frac{k+m_2+l-1}{l}}
\right)ds\\
&\le \frac{k(k-1)\chi l}{(k+m_2-1)(k+m_2+l-1)} C_P^{\frac{l}{k+m_2+l-1}}
K^{-\frac{k+m_2-1}{l}}
\left[\lVert v_0\rVert^q_{1,1-\frac{1}{q}} +\int_0^t e^s \left(\int_\Omega |f(\cdot,s)|^{\frac{k+m_2+l-1}{l}}\right)ds\right]\\
&\leq \const{G2} + \frac{k(k-1)\chi l K}{(k+m_2-1)(k+m_2+l-1)} C_P^{\frac{l}{k+m_2+l-1}} \int_0^t e^s \left(\int_\Omega (u+1)^{k+m_2+l-1}\right)ds \quad \text{on } (0,\TM).
\end{split}
\end{equation}
On the other hand, by virtue of \eqref{Y2} and \eqref{ParabSt1}, estimate \eqref{Stima1} can be rewritten as 
\begin{equation}\label{ParabSt3}
\begin{split}
y'(t)&\leq -\frac{4k(k-1)}{(k+m_1-1)^2} \int_\Omega |\nabla (u+1)^{\frac{k+m_1-1}{2}}|^2 + \frac{k(k-1)\chi K C_P^{\frac{l}{k+m_2+l-1}}}{k+m_2+l-1}
\int_\Omega (u+1)^{k+m_2+l-1}\\ 
& \quad +  \frac{k(k-1)\chi l}{(k+m_2-1)(k+m_2+l-1)}  \left(K C_P^{\frac{l}{k+m_2+l-1}}\right)^{-\frac{k+m_2-1}{l}} \int_\Omega |\Delta v|^{\frac{k+m_2+l-1}{l}}\\ 
&\quad +\left(\epsilon-\frac{kb}{2^{\beta-1}}\right) \int_\Omega (u+1)^{k+\beta-1}
- \frac{kc}{2^{\delta-1}} \int_\Omega (u+1)^{k-1} \int_\Omega (u+1)^{\delta} + \const{S1} \quad \mbox{on}\ (0, \TM).
\end{split}
\end{equation}
Now we add to both sides of \eqref{ParabSt3} the term $y(t)$, and successively we multiply by $e^t$. Since $e^t y'(t) + e^t y(t) =\frac{d}{dt}(e^t y(t))$, an integration over $(0,t)$ provides for all $t \in (0,\TM)$
\begin{equation*}
\begin{split}
e^t y(t) \leq y(0) +\int_0^t e^s &\left(-\frac{4k(k-1)}{(k+m_1-1)^2}\int_\Omega |\nabla (u+1)^{\frac{k+m_1-1}{2}}|^2 
+ \frac{k(k-1)\chi K C_P^{\frac{l}{k+m_2+l-1}}}{k+m_2+l-1} \int_\Omega (u+1)^{k+m_2+l-1} \right.\\
& \quad \left. + \frac{k(k-1)\chi l}{(k+m_2-1)(k+m_2+l-1)}  \left(K C_P^{\frac{l}{k+m_2+l-1}}\right)^{-\frac{k+m_2-1}{l}} \int_\Omega |\Delta v|^{\frac{k+m_2+l-1}{l}}\right.\\
&\quad \left.  +\left(\epsilon-\frac{kb}{2^{\beta-1}}\right)\int_\Omega (u+1)^{k+\beta-1}+\int_\Omega (u+1)^k
- \frac{kc}{2^{\delta-1}} \int_\Omega (u+1)^{k-1} \int_\Omega (u+1)^{\delta} + \const{S1}\right)\, ds, 
\end{split}
\end{equation*}
which, taking into account bound \eqref{ParabSt2} and applying Young's inequality  to $\int_\Omega (u+1)^k$ to control $\int_\Omega (u+1)^{k+\beta -1}$ (recall $\beta>1$), gives 
\begin{equation}\label{ParabSt4}
\begin{split}
e^t y(t) \leq y(0) +\int_0^t e^s &\left(-\frac{4k(k-1)}{(k+m_1-1)^2}\int_\Omega |\nabla (u+1)^{\frac{k+m_1-1}{2}}|^2 \right.\\
& \quad \left. +\frac{k(k-1)\chi K}{k+m_2-1} C_P^{\frac{l}{k+m_2+l-1}} \int_\Omega (u+1)^{k+m_2+l-1}  +\left(\epsilon-\frac{kb}{2^{\beta-1}}\right)  \int_\Omega (u+1)^{k+\beta-1}\right.\\ 
& \quad \left.- \frac{kc}{2^{\delta-1}} \int_\Omega (u+1)^{k-1} \int_\Omega (u+1)^{\delta} + \const{sf1}\right)\, ds \quad \textrm{on } (0,\TM).
\end{split}
\end{equation}
We now rely again on estimate \eqref{Interpolation1} for $\hat{c}=\frac{k(k-1)\chi K}{k+m_2-1} C_P^{\frac{l}{k+m_2+l-1}}$; by following similar arguments as the case $\tau=0$ above, we have by using the first bound in  \eqref{Interpolation1}
\begin{equation}\label{StimaFinale} 
e^t \int_\Omega (u+1)^k\leq \const{sd} + \const{a_6e}(e^t -1) \quad \textrm{for all }  t \in (0,\TM),    
 \end{equation}
so concluding. The second bound of \eqref{Interpolation1} is, indeed, invoked to establish the limit case and in this situation the largeness hypothesis on $c$ has to be as well taken in mind.
% with
%\begin{equation*}
%\int_\Omega u^k \leq \const{sdf}\quad \textrm{on } (0,\TM).
%\end{equation*}
\end{proof}
\end{lemma}
\subsection{A priori boundedness of $u$ for model \eqref{problem_grad}}
\begin{lemma}\label{StimaGeneraleP}
Let the assumptions of Lemma \ref{ThetaSigma} be valid and let $k_*>1$ be 
the value described in Convention \ref{Convenzione}. Then, for $k>k_*$, the functional $y(t)$ defined in \eqref{Func} fulfills for every $\epsilon>0$ and for $F_{m_2}(u)$ as in \eqref{Def_Fm2}
\begin{equation}\label{Stima1P}
\begin{split}
y'(t)&\leq -\frac{4k(k-1)}{(k+m_1-1)^2} \int_\Omega |\nabla (u+1)^{\frac{k+m_1-1}{2}}|^2 - k(k-1) \chi \int_\Omega F_{m_2}(u)\Delta v \\
&\quad+ k (a +b|\Omega|+\epsilon) \int_\Omega (u+1)^{k+\alpha-1} +\left(\epsilon- \frac{kb}{2^{\alpha+\beta-2}}\right) \int_\Omega (u+1)^{k+\alpha-1} \int_\Omega (u+1)^{\beta}\\
&\quad - kc \left(\frac{k-1+\delta}{\delta}\right)^{-\delta} \int_\Omega  
|\nabla (u+1)^{\frac{k-1+\delta}{\delta}}|^{\delta} + \const{Aa} 
\quad \mbox{on}\ (0, \TM).
\end{split}
\end{equation}
\begin{proof}
Retracing the reasoning in Lemma \ref{StimaGenerale}, from the first equation in \eqref{problem_grad} we deduce for all $k>k_*$ and $t \in (0, \TM)$ that 
\begin{equation*}
\begin{split}
y'(t) &= -\frac{4k(k-1)}{(k+m_1-1)^2} \int_\Omega |\nabla (u+1)^{\frac{k+m_1-1}{2}}|^2 - k(k-1) \chi \int_\Omega F_{m_2}(u)\Delta v\\ 
&\quad +ka \int_\Omega u^\alpha (u+1)^{k-1} - kb  \int_\Omega u^{\alpha} (u+1)^{k-1}  \int_\Omega u^{\beta} - kc \left(\frac{k-1+\delta}{\delta}\right)^{-\delta} \int_\Omega  
|\nabla (u+1)^{\frac{k-1+\delta}{\delta}}|^{\delta}.
\end{split}
\end{equation*}
From $u \leq u+1$ and inequality \eqref{ABIn}, the previous estimate becomes on $(0,\TM)$
\begin{equation}\label{Stima3P}
\begin{split}
y'(t) &\leq -\frac{4k(k-1)}{(k+m_1-1)^2} \int_\Omega |\nabla (u+1)^{\frac{k+m_1-1}{2}}|^2 - k(k-1) \chi \int_\Omega F_{m_2}(u)\Delta v\\ 
&\quad +k (a +b|\Omega|) \int_\Omega (u+1)^{k+\alpha-1} 
- \frac{kb}{2^{\alpha+\beta-2}}  \int_\Omega (u+1)^{k+\alpha-1}  \int_\Omega (u+1)^{\beta} \\
& \quad + \frac{kb}{2^{\beta-1}} \int_\Omega (u+1)^{k-1}  \int_\Omega (u+1)^{\beta}
- kc \left(\frac{k-1+\delta}{\delta}\right)^{-\delta} \int_\Omega  
|\nabla (u+1)^{\frac{k-1+\delta}{\delta}}|^{\delta}.
\end{split}
\end{equation}
Let us deal with the second last integral of the right-hand side of \eqref{Stima3P}; we can apply the Young inequality to write 
\[
\int_\Omega (u+1)^{\beta}\leq \epsilon \int_\Omega (u+1)^{k+\alpha-1}+\const{TT} \quad \textrm{on } t\in (0,\TM),
\]
and by virtue of $\alpha\geq 1$ 
\begin{equation}\label{YP1}
\frac{kb}{2^{\beta-1}} \int_\Omega (u+1)^{k-1}  \int_\Omega (u+1)^{\beta} \leq \epsilon\int_\Omega (u+1)^{k+\alpha-1}\int_\Omega (u+1)^{\beta}+ \epsilon \int_\Omega (u+1)^{k+\alpha-1}+\const{axy}   \quad \textrm{on } (0,\TM).
\end{equation}
By putting \eqref{YP1} into estimate \eqref{Stima3P}, we conclude the proof. 
\end{proof}
\end{lemma}
\begin{lemma}\label{LemmaProblemGradEll2}
Let the assumptions of Lemma \ref{ThetaSigma} be valid, let $k_*>1$ be 
the value described in Convention \ref{Convenzione} and let $k>k_*$. Additionally, let $f$ comply with \eqref{Def_f} and let $\tau\in \{0,1\}$, $c\geq 0$, $\alpha\geq 1$, $\beta \geq \max\{1,m_1-\alpha\}$ and  $1\leq\delta \leq 2$ satisfy 
\begin{itemize}
\item[$\triangleright$] if $c=0$, the condition
\begin{equation}\label{czeroG}
 \beta> \max\left\{\frac{n(m_2-m_1+l)+2(m_2+l-\alpha)}{2}, \frac{n(\alpha-m_1)}{2}\right\};
\end{equation}
\item[$\triangleright$] if $c>0$, either the constrains
\begin{equation}\label{cdivzeroG1}
 \beta>  \frac{n(\alpha-m_1)}{2} \quad \textrm{and} \quad \delta>\frac{n(m_2+l)}{n+1};
\end{equation}
%\item if $c>0$ the boundsç
or
\begin{equation}\label{cdivzeroG12}
 \beta>  \frac{n(m_2-m_1+l)+2(m_2+l-\alpha)}{2} \quad \textrm{and} \quad \delta> \frac{n \alpha}{n+1}.
\end{equation}
\end{itemize}
Then $u \in L^{\infty}((0,\TM), L^k(\Omega))$ for all $k>k_*$.
%%%% CASI LIMITI

Further, let $H_{1,\alpha}=1$, for $\alpha=1$, and $H_{1,\alpha}=0$, for $\alpha\neq 1$.
The same conclusion continues to be valid 
\begin{itemize}
\item[$\triangleright$] either for $c=0$ in these circumstances
\begin{equation}\label{czeroGLimitCase}
\begin{dcases}
\displaystyle \beta= \frac{n(m_2-m_1+l)+2(m_2+l-\alpha)}{2} > \frac{n(\alpha-m_1)}{2} \\ \quad  \quad \quad \quad  \quad \quad \quad  \quad \quad \quad  \quad \quad \textrm{if } \displaystyle b>  \frac{2^{\alpha+\beta-2}(k-1)\chi K C_0}{k+m_2-1}\left(\tau C_P^\frac{l}{k+m_2+l-1}+(1-\tau)\right), \\
\displaystyle \beta=  \frac{n(\alpha-m_1)}{2} > \frac{n(m_2-m_1+l)+2(m_2+l-\alpha)}{2}\\ \quad  \quad \quad \quad  \quad 
 \textrm{if }  \displaystyle
a<\frac{b}{\mathcal{C}_0 2^{\alpha +\beta-2}}-\frac{\tau H_{1,\alpha}}{k}-b|\Omega|
\quad \textrm{provided} \quad \frac{b}{\mathcal{C}_0 2^{\alpha +\beta-2}}-\frac{\tau H_{1,\alpha}}{k}-b|\Omega|>0, \\
\displaystyle \beta=\frac{n(m_2-m_1+l)+2(m_2+l-\alpha)}{2}= \frac{n(\alpha-m_1)}{2}
\\ \quad  \quad \quad \quad  \quad 
\textrm{if }  \displaystyle
b>2^{\alpha+\beta-2}\mathcal{C}_0\left[\frac{(k-1)\chi K}{k+m_2-1}\left(\tau \mathcal{C}_P^\frac{l}{k+m_2+l-1}+(1-\tau)\right)+a+b|\Omega|+\frac{\tau H_{1,\alpha}}{k}\right],
\end{dcases}
\end{equation}
\item[$\triangleright$] or for
\begin{equation}\label{cdivzeroG1LimitCase}
\begin{dcases}
 \displaystyle \beta>  \frac{n(\alpha-m_1)}{2} \textrm{ and } \delta=\frac{n(m_2+l)}{n+1} \\ \quad \quad \quad  \quad \quad  \textrm{if } \displaystyle c\geq \frac{(k-1)\chi K}{k+m_2-1}(2 C_{GN})^{\sigma_3} (M_1+|\Omega|)^{\frac{m_2+l}{n+1}}\left(\frac{k-1+\delta}{\delta}\right)^{\delta}\left(\tau C_P^\frac{l}{k+m_2+l-1}+(1-\tau)\right),\\ 
\displaystyle \beta>  \frac{n(m_2-m_1+l)+2(m_2+l-\alpha)}{2} \textrm{ and } \delta= \frac{n \alpha}{n+1} \\
\quad \quad \quad \quad \quad   \textrm{if }\displaystyle 
c > (a+b |\Omega|)(2 C_{GN})^{\sigma_2} (M_1+|\Omega|)^{\frac{\alpha}{n+1}}\left(\frac{k-1+\delta}{\delta}\right)^{\delta}.
\end{dcases}
\end{equation}
\end{itemize}
\begin{proof}
As before we distinguish the cases $\tau=0$ and $\tau=1$, and moreover we give the details only for one limit case; the reasoning for the other situations are adaptations of what already done for problem \eqref{problem_int}, and precisely in Lemma \ref{GNandInterp} and Lemma \ref{GNandInterpBisLemma}.
\subsubsection*{$\triangleright$ $\tau=0$}
By inserting estimate \eqref{StimaFm2} in \eqref{Stima1P}, we derive the following bound
\begin{equation}\label{Stima2P}
\begin{split}
y'(t)&\leq -\frac{4k(k-1)}{(k+m_1-1)^2} \int_\Omega |\nabla (u+1)^{\frac{k+m_1-1}{2}}|^2 + \frac{k(k-1)\chi K}{k+m_2-1} \int_\Omega (u+1)^{k+m_2+l-1}\\
&\quad+ k \left(a +b|\Omega|+\epsilon\right) \int_\Omega (u+1)^{k+\alpha-1} +\left(\epsilon- \frac{kb}{2^{\alpha+\beta-2}}\right) \int_\Omega (u+1)^{k+\alpha-1} \int_\Omega (u+1)^{\beta}\\
&\quad - kc \left(\frac{k-1+\delta}{\delta}\right)^{-\delta} \int_\Omega  
|\nabla (u+1)^{\frac{k-1+\delta}{\delta}}|^{\delta} + \const{Aa} 
\quad \mbox{on}\ (0, \TM).
\end{split}
\end{equation}
For $c=0$, the assumption on $\beta$ in \eqref{czeroG} allows to exploit the first inequalities in \eqref{InterpolationP1} and in \eqref{InterpolationP2} in order to estimate the second and the third integrals on the right-hand side of \eqref{Stima2P}, so arriving at 
\begin{equation}\label{ES}
y'(t) \leq - \const{Hh} \int_\Omega |\nabla (u+1)^{\frac{k+m_1-1}{2}}|^2 + \const{Ssi}  \quad \mbox{on}\ (0, \TM).
\end{equation}
When $c>0$, by starting from inequality \eqref{Stima2P} and by taking into account the restrictions on $\beta$ and $\delta$ in \eqref{cdivzeroG1}, the integrals 
$\int_\Omega (u+1)^{k+\alpha-1}$ and 
$\int_\Omega (u+1)^{k+m_2+l-1}$ can be estimated by using the first bounds in \eqref{InterpolationP2} and in \eqref{GN4}, which, plugged into \eqref{Stima2P}, give an estimate similar to \eqref{ES}. Alternatively, by virtue of \eqref{cdivzeroG12} we can also apply the first inequalities in \eqref{GN3} and in \eqref{InterpolationP1} in order to control the integrals $\int_\Omega (u+1)^{k+\alpha-1}$ and 
$\int_\Omega (u+1)^{k+m_2+l-1}$, leading to bound  \eqref{ES} up to constants. With these preparations, we conclude by reasoning as in Lemma \ref{LemmaProblemIntEll1} for the case $\tau=0$.

Let us now deal with the first limit case in \eqref{czeroGLimitCase}: 
%%We insert inequality \eqref{StimaFm2} into bound \eqref{Stima3P}, and this yields
%\begin{equation*}
%\begin{split}
%y'(t) &\leq -\frac{4k(k-1)}{(k+m_1-1)^2} \int_\Omega |\nabla (u+1)^{\frac{k+m_1-1}{2}}|^2 + \frac{k(k-1)\chi K}{k+m_2-1} \int_\Omega (u+1)^{k+m_2+l-1}\\ 
%&\quad +k (a +b|\Omega|) \int_\Omega (u+1)^{k+\alpha-1} 
%- \frac{kb}{2^{\alpha+\beta-2}}  \int_\Omega (u+1)^{k+\alpha-1}  \int_\Omega (u+1)^{\beta} \\
%& \quad + \frac{kb}{2^{\beta-1}} \int_\Omega (u+1)^{k-1}  \int_\Omega (u+1)^{\beta}
%- kc \left(\frac{k-1+\delta}{\delta}\right)^{-\delta} \int_\Omega  
%|\nabla (u+1)^{\frac{k-1+\delta}{\delta}}|^{\delta}\\
%&\leq -\frac{4k(k-1)}{(k+m_1-1)^2} \int_\Omega |\nabla (u+1)^{\frac{k+m_1-1}{2}}|^2 + \frac{k(k-1)\chi K}{k+m_2-1} \int_\Omega (u+1)^{k+m_2+l-1}\\ 
%&\quad +k (a +b|\Omega|) \int_\Omega (u+1)^{k+\alpha-1} 
%+\left(\frac{\epsilon}{2}- \frac{kb}{2^{\alpha+\beta-2}} \right) \int_\Omega (u+1)^{k+\alpha-1}  \int_\Omega (u+1)^{\beta} \\
%& \quad - kc \left(\frac{k-1+\delta}{\delta}\right)^{-\delta} \int_\Omega  
%|\nabla (u+1)^{\frac{k-1+\delta}{\delta}}|^{\delta} \quad \textrm{for all } t \in (0,\TM).
%\end{split}
%\end{equation*}
since $\beta>\frac{n}{2}(\alpha-m_1)$, by plugging the second estimate in \eqref{InterpolationP1} with $\hat{c}=\frac{k(k-1)\chi K C_0}{k+m_2-1}$ and the first in \eqref{InterpolationP2} with $\hat{c}=k(a+b|\Omega|+\epsilon)$ into bound \eqref{Stima2P}, we obtain
\begin{equation}\label{Stima2Pp}
\begin{split}
y'(t)&\leq -\frac{2k(k-1)}{(k+m_1-1)^2} \int_\Omega |\nabla (u+1)^{\frac{k+m_1-1}{2}}|^2 \\
&\qquad +\left(\epsilon+\frac{k(k-1)\chi K C_0}{k+m_2-1}- \frac{kb}{2^{\alpha+\beta-2}}\right) \int_\Omega (u+1)^{k+\alpha-1} \int_\Omega (u+1)^{\beta}+ \const{Aaa} 
\quad \mbox{on}\ (0, \TM).
\end{split}
\end{equation}
Since $b > \frac{2^{\alpha+\beta-2}(k-1)\chi K C_0}{k+m_2-1}$,  inequality \eqref{Stima2Pp} gives an estimate similar to \eqref{ES}, therefore the conclusion is obtained. 
\subsubsection*{$\triangleright$ $\tau=1$}
%\end{proof}
%\end{lemma}
%\begin{lemma}\label{LemmaProblemGradParab2}
%Let the assumptions of Lemma \ref{ThetaSigma} be valid and let $k_*>1$ be as in Convention \ref{Convenzione} the value therein described. Additionally, let $\tau=1$. Assume that $\alpha, \beta \geq 1$ and $1\leq\delta \leq 2$ satisfy conditions 
%\eqref{czeroG} if $c=0$ and either \eqref{cdivzeroG1} or \eqref{cdivzeroG12} if  $c>0$, respectively.
%Then $u \in L^{\infty}((0,\TM), L^k(\Omega))$ for all $k>k_*$.
%\begin{proof}
By means of \eqref{ParabSt1} and \eqref{ParabSt2}, after some adjustments 
we deduce the following estimate
\begin{equation*}\label{PPGrad}
\begin{split}
e^t y(t) \leq y(0) +\int_0^t e^s &\left(-\frac{4k(k-1)}{(k+m_1-1)^2}\int_\Omega |\nabla (u+1)^{\frac{k+m_1-1}{2}}|^2 \right.\\
& \quad \left.+\frac{k(k-1)\chi K}{k+m_2-1} C_P^{\frac{l}{k+m_2+l-1}} \int_\Omega (u+1)^{k+m_2+l-1} \right.\\
&\quad \left. + k \left(a +b|\Omega|+\epsilon\right) \int_\Omega (u+1)^{k+\alpha-1} +\left(\epsilon- \frac{kb}{2^{\alpha+\beta-2}}\right) \int_\Omega (u+1)^{k+\alpha-1} \int_\Omega (u+1)^{\beta}\right.\\
& \quad \left.- kc \left(\frac{k-1+\delta}{\delta}\right)^{-\delta} \int_\Omega  
|\nabla (u+1)^{\frac{k-1+\delta}{\delta}}|^{\delta} +\int_\Omega (u+1)^k+ \const{AAa} \right)\, ds \quad \textrm{on } (0,\TM).
\end{split}
\end{equation*}
To treat the term $\int_\Omega (u+1)^k$, the analysis is more manageable for $\alpha>1$, but less straightforward for $\alpha=1$. Indeed, for $\alpha>1$, the Young inequality allows to write $\int_\Omega (u+1)^{k}\leq \epsilon\int_\Omega (u+1)^{k+\alpha-1}+\const{pdd}$. If, conversely, $\alpha=1$ the presence of the nonlinearity $\int_\Omega  
|\nabla (u+1)^{\frac{k-1+\delta}{\delta}}|^{\delta}$, is crucial: in fact, for $c>0$ the gradient term can be used to control $\int_\Omega (u+1)^k$, exactly by exploiting the first in \eqref{GN3} for $\alpha=1$ (note $\delta\geq1>\frac{n}{n+1}$). For $c=0$ (and $\alpha=1$), we make use of \eqref{InterpolationP2} and the resulting term $\left(k \left(a +b|\Omega|+\epsilon\right)+1\right)\int_\Omega (u+1)^k$ can be absorbed by the remaining ones. 
We establish the claim by exactly retracing what done some lines above 
for the cases $c=0$ and $c>0$ and we conclude by following similar argument to Lemma \ref{LemmaProblemIntEll1} for the case $\tau=1$.

%\begin{equation}\label{StimaFinale} 
%e^t \int_\Omega (u+1)^k\leq \const{Ssd} + \const{A_6e}(e^t -1) \quad \textrm{for all }  t \in (0,\TM),    
 %\end{equation}
%henceforth, we conclude that
%\begin{equation*}
%\int_\Omega u^k \leq \const{svdf}\quad \textrm{on } (0,\TM).
%\end{equation*}
%\end{proof}
%\end{lemma}
%\section{The limit cases for problems \eqref{problem_int} and \eqref{problem_grad}}\label{LimitCases}
%\subsection{Analysis of problem \eqref{problem_grad}}
%\begin{lemma}\label{PreCorGrad1}
%Let the assumptions of Lemma \ref{ThetaSigma} be valid and let $k_*>1$ be as in Convention %\ref{Convenzione} the value therein described, $c=0$ and 
%$\beta=\frac{n(m_2-m_1+l)+2(m_2+l-\alpha)}{2}>\max\left\{m_1-\alpha,\frac{n(\alpha-m_1)}{2}\right\}$.
%Then the same claims of Lemma \ref{LemmaProblemGradEll2} and Lemma \ref{LemmaProblemGradParab2} hold %provided that
%\begin{itemize}
%\item $b > \frac{2^{\alpha+\beta-2}(k-1)\chi K C_0}{k+m_2-1} C_P^{\frac{l}{k+m_2+l-1}}$ if $\tau=1$,
%\end{itemize}
%where $C_0, C_P$ are given \textcolor{blue}{PRECISA...ESPRESS ENORME, DIVERSA DA CASO A CASO} in Lemma \ref{GagliardoIneqLemma} and Lemma \ref{lem:MaxReg}, respectively.
%\begin{proof}
\end{proof}
\end{lemma}
We can now conclude and give the proofs of the claimed results.
\subsection{Proofs of Theorems \ref{MainTheoremInt} and \ref{MainTheoremGrad}, and Corollaries \ref{corollaryToTheo1} and \ref{corollaryToTheo2} }\label{SectionProofTheorem}
With the aid of Convention \ref{Convenzione} let $k_*>1$ be the value in Lemma \ref{ThetaSigma}.
Additionally, accordingly to the first item in Remark \ref{RemarkCostanti}, let $\mathcal{C}_{GN}=\mathcal{C}_{GN}(k_*), \mathcal{C}_0=\mathcal{C}_0(k_*)$ and $\mathcal{C}_P=\mathcal{C}_P(k_*)$. Since the values $\mathcal{C}_{GN}=\mathcal{C}_{GN}(k), \mathcal{C}_0=\mathcal{C}_0(k)$ and $\mathcal{C}_P=\mathcal{C}_P(k)$ continuously  depend on $k$, restrictions \eqref{LimitCaseCInt} of Corollary \ref{corollaryToTheo1} and \eqref{LimitA}--\eqref{LimitD} of Corollary \ref{corollaryToTheo2} are such that relations \eqref{LimitInt1} in Lemma \ref{LemmaProblemIntEll1} and
\eqref{czeroGLimitCase} and \eqref{cdivzeroG1LimitCase} in Lemma \ref{LemmaProblemGradEll2} are satisfied  for some $k>k_*$. Thereafter, $u \in L^{\infty}((0,\TM);L^k(\Omega))$ and subsequently 
elliptic and parabolic regularity results, in conjunction with Sobolev embeddings, imply $\nabla v\in L^\infty((0,\TM);L^\infty(\Omega))$. In addition, the sources in problems \eqref{problem_int} and \eqref{problem_grad} belong to $L^\infty((0,\TM);L^{k}(\Omega))$. 
As a consequence, from \cite[Lemma A.1]{TaoWinkParaPara} we have 
$u \in L^{\infty}((0,\TM);L^{\infty}(\Omega))$, so concluding by applying Lemma \ref{ExtensionLemma}. Theorems \ref{MainTheoremInt} and \ref{MainTheoremGrad} similarly follow taking in mind their hypotheses, and lemmata \ref{LemmaProblemIntEll1}, \ref{LemmaProblemGradEll2} and \ref{ExtensionLemma}.

\qed

\appendix
\section{Some other scenarios toward boundedness}
In this appendix we briefly dedicate ourselves to analyze further situations where boundedness of solutions to models \eqref{problem_int} and \eqref{problem_grad} is ensured.
\subsection{Situations for model \eqref{problem_int}}\label{RemarkLinear}
\begin{enumerate}
\item \emph{Absence of nonlocal term.} When $\beta>m_2+l$, in order to achieve the crucial inequality \eqref{uk_estimate4} (if $\tau=0$), or \eqref{StimaFinale} (if $\tau=1$), it is sufficient to control the other terms in \eqref{Sf1} or \eqref{ParabSt4} (depending on the value of $\tau$) by means of $\int_\Omega (u+1)^{k+\beta-1}$, coming from the logistic. In the specific for $\beta=m_2+l$, the achievement of \eqref{uk_estimate4} (or \eqref{StimaFinale}) is not direct and the extra term
\begin{equation*}
\left(\frac{k(k-1)\chi K}{k+m_2-1}\left( \tau C_P^{\frac{l}{k+\beta-1}}+(1-\tau)\right)+\epsilon
-\frac{kb}{2^{\beta-1}} \right) \int_\Omega (u+1)^{k+\beta-1}
\end{equation*} 
would take part to the computations. Subsequently, restriction in \eqref{RestricionOnbModelNonlocal} is required.

As to the linear case $m_1=m_2=l=1$, the study is developed with the simplified functional $y(t)=\int_\Omega u^k$, for both $\tau=0$ and $\tau=1$. 
Analogous reasonings show that expression 
\begin{equation*}
\left(\frac{k(k-1)\chi K}{k+m_2-1}\left( \tau C_P^{\frac{l}{k+\beta-1}}+(1-\tau)\right)+\epsilon-\frac{kb}{2^{\beta-1}}\right) \int_\Omega (u+1)^{k+\beta-1}
\end{equation*} 
actually reduces to
\begin{equation*}
\left((k-1)\chi K\left( \tau C_P^{\frac{1}{k+\beta-1}}+(1-\tau)\right)+\epsilon-kb\right) \int_\Omega u^{k+\beta-1}.
\end{equation*}
Recalling Remark \ref{MainTheorem-CorollaryInt}, relation \eqref{RestricionOnbModelNonlocalLinear} is established by imposing the above expression in brackets negative, 
for some $k>\frac{n}{2}$ and $\beta=2$.

\item \emph{Absence of logistic.} Further conditions are either
\[
m_2+l<m_1+ \frac{2}{n},
\]
or if $m_2+l=m_1+ \frac{2}{n}$
\[
\chi < \frac{4(k_*+m_2-1)}{K (k_*+m_1-1)^2} (2 \mathcal{C}_{GN})^{-\frac{2(k_*+m_1+\frac{2}{n}-1)}{k_*+m_1-1}} (M_0 + |\Omega|)^{-\frac{2}{n}} \left(\tau \mathcal{C}_P^{-\frac{l}{k_*+m_1+\frac{2}{n}-1}}+(1-\tau)\right),
\]
where $M_0$ is given in bound \eqref{MassBounded_int}.
\end{enumerate}

\subsection{Situations for model \eqref{problem_grad}}\label{RemarkGrad}
\begin{enumerate}
\item \emph{Absence of the gradient term.} For $\tau=0$ in problem \eqref{problem_grad}, the conditions below are as well sufficient to have boundedness of related solutions:
\[
\alpha, \beta \geq 1 \quad \text{and} \quad m_1>\max\left\{\alpha-\frac{2}{n}, m_2+l-\frac{2}{n}\right\}.
\]
The restriction $m_1>m_2+l-\frac{2}{n}$ has been already discussed in $\S$\ref{RemarkLinear}, while 
$m_1>\alpha-\frac{2}{n}$ is required to apply the Gagliardo--Nirenberg and Young's inequalities to derive this estimate for $\hat{c}>0$
\begin{equation*}
\hat{c} \int_\Omega (u+1)^{k+\alpha-1} \leq \frac{k(k-1)}{(k+m_1-1)^2} \int_\Omega 
 |\nabla (u+1)^{\frac{k+m_1-1}{2}}|^2 + \const{SFsf}  \quad \text{on }(0,\TM),
\end{equation*}
successively employed in \eqref{Stima2P}. (The fully parabolic case is similar.)
In particular, we underline that $m_1>m_2+l-\frac{2}{n}$ improves the related assumption in \cite[Theorem 1.2]{TANAKAEtAl2021} for $l \in \left(0,\frac{1}{n}\right]$.

\item \emph{Presence of the gradient term.} Boundedness is also obtained for:  
\begin{itemize}
\item[$\triangleright$] $\alpha, \beta \geq 1$ and $\max\left\{\frac{n\alpha}{n+1}, \frac{n(m_2+l)}{n+1}\right\}< \delta \leq 2$. 
\item[$\triangleright$] $1 \leq \alpha < m_1 + \frac{2}{n}$, $\beta \geq 1$ and $\frac{n(m_2+l)}{n+1}< \delta \leq 2$. 
\end{itemize}
\end{enumerate}

\subsubsection*{\bf\textit{\quad Acknowledgments}}
The authors are members of the Gruppo Nazionale per l'Analisi Matematica, la Probabilit\`a e le loro Applicazioni (GNAMPA) of the Istituto Nazionale di Alta Matematica (INdAM) and they are partially supported by the research project {\em  Analysis of PDEs in connection with real phenomena} (2021, Grant Number: F73C22001130007), funded by \href{https://www.fondazionedisardegna.it/}{Fondazione di Sardegna}.
SF and GV are also partially supported by MIUR (Italian Ministry of Education, University and Research) Prin 2022 \emph{Nonlinear differential problems with applications to real phenomena} (Grant Number: 2022ZXZTN2).
SF also acknowledges financial support by INdAM-GNAMPA project \emph{Problemi non lineari di tipo stazionario ed evolutivo} (CUP E53C23001670001). GV is also partially supported by INdAM-GNAMPA project \emph{Equazioni differenziali alle derivate parziali nella modellizzazione di fenomeni reali} (CUP E53C22001930001) and \emph{Studio di modelli nelle scienze della vita} (CUP J55F21004240001, ``DM 737/2021 risorse 2022--2023'', founded by European Union--NextGenerationEU). 
RDF acknowledges financial support by PNRR e.INS Ecosystem of Innovation for Next Generation Sardinia (CUP F53C22000430001, codice MUR ECS0000038). 

\subsubsection*{\bf\textit{\quad Conflict of interests}}
The authors declare that they have no conflict of interests.

%\bibliography{Bibliography}{}
%\bibliographystyle{abbrv}
\end{document}